\documentclass{amsart}

\usepackage[all]{xy}
 \usepackage[usenames,dvipsnames]{pstricks}
\usepackage{ebezier}
\usepackage[pagewise]{lineno}

\usepackage{latexsym}\usepackage{amsfonts}\usepackage{amssymb}\usepackage{amsthm}
\usepackage{amsmath}\usepackage{newlfont}\usepackage{graphicx}\usepackage{doc}
\usepackage{color}\usepackage{multicol}\usepackage{epic}\usepackage{eepic}
\usepackage{ebezier}
\newtheorem{thrm}{Theorem}
\newtheorem{lem}[thrm]{Lemma}
\newtheorem{cor}[thrm]{Corollary}
\newtheorem{defn}[thrm]{Definition}
\newtheorem{rem}[thrm]{Remark}
\newtheorem{prop}[thrm]{Proposition}
\include{epsf}

\def \Dj{\mbox{\raise0.3ex\hbox{-}\kern-0.4em D}}

\DeclareMathOperator{\limarr}{\displaystyle \lim_{\longrightarrow}}
\newcommand{\Id}{\operatorname{Id}}
\newcommand{\IM}{\operatorname{Im}}
\newcommand{\Crit}{\operatorname{Crit}}
\newcommand{\sing}{\operatorname{sing}}

\newcommand{\PSS}{\operatorname{PSS}}
\newcommand{\PD}{\operatorname{PD}}
\newcommand{\ev}{\operatorname{ev}}

\begin{document}

\title{Spectral Invariants in Lagrangian Floer homology of open subset}

\author{Jelena Kati\'c}
\author{Darko Milinkovi\'c} 
\author{Jovana Nikoli\'c}

\address{Matemati\v{c}ki fakultet, Studentski trg 16, 11000
Belgrade, Serbia}

\email{jelenak@matf.bg.ac.rs, milinko@matf.bg.ac.rs,jovanadj@matf.bg.ac.rs}

\thanks{This work is
partially supported by Ministry of Education and Science
of Republic of Serbia Project \#ON174034.}

\begin{abstract}
We define and investigate spectral invariants for Floer homology $HF(H,U:M)$ of an open subset $U\subset M$ in $T^*M$, defined by Kasturirangan and Oh as a direct limit of Floer homologies of approximations. We define a module structure product on $HF(H,U:M)$ and prove the triangle inequality for invariants with respect to this product. We also prove the continuity of these invariants and compare them with spectral invariants for periodic orbits case in $T^*M$.
\end{abstract}

\maketitle


MSC[2010] Primary 53D12, Secondary 53D40


\section{Introduction}

\subsection{Spectral invariants in cotangent bundles}
Spectral invariants in Symplectic Topology in terms of generating functions for Lagrangian submanifolds of cotangent bundles were introduced by Viterbo in \cite{V}. If $E\rightarrow M$ is a smooth vector bundle over a compact smooth manifold $M$, $S:E\rightarrow\mathbb{R}$ a generic smooth function and
$$
\Sigma_S:=\{e\in E \mid d_{vert}S(e)=0\}
$$
(here $d_{vert}S$ denotes the derivative along the fibre), then
$$
i_S:\Sigma_S\rightarrow T^*M, \qquad i_S(e):=dS(e) 
$$
is a smooth Lagrangian immersion. It is known that all Hamiltonian deformations of zero section can be generated by some function $S$ in this way~\cite{LS,C1,C2}. Viterbo defined spectral invariants as a certain minimax values of $S$. He used them to prove several important results about Hamiltonian diffeomorphisms.

In~\cite{O1,O2} Oh defined spectral invariants for the case of cotangent bundle using the ``homologically visible" critical values of the action functional
$$a_H(x):=\int_x\theta-\int_0^1H(x(t),t)dt,$$ 
where $\theta$ is the Liouville $1-$form on $T^*M$. More precisely, let $O_M$ be a zero section of $T^*M$ and
$L=\phi^1_H(O_M)$, where $\phi^1_H$ is a time--one--map generated by a Hamiltonian $H$. Let
$HF_*^{\lambda}(O_M,\phi^1_H(O_M))$ denotes the filtrated homology defined via filtrated Floer complex:
$$CF_*^{\lambda}(O_M,\phi^1_H(O_M)):=\mathbb Z_2\langle\{x\in\Crit(a_H)\mid a_H(x)<\lambda\}\rangle.$$
These homology groups are well defined since the boundary map preserves the filtration:
$$\partial:CF_*^{\lambda}(O_M,\phi^1_H(O_M))\to CF_*^{\lambda}(O_M,\phi^1_H(O_M)),$$
due to well defined action functional that decreases along its ``negative gradient flows".
For a singular homology class $\alpha\in H_*(M,\mathbb{Z}_2)$ define
$$\sigma(\alpha,H):=
\inf\{{\lambda}\in\mathbb R\mid F_H(\alpha)\in\IM(\imath_*^{\lambda})\}$$ where
$$\imath_*^{\lambda}:HF_*^{\lambda}(O_M,\phi^1_H(O_M))
\to HF_*(O_M,\phi^1_H(O_M))$$ is the homomorphism induced by inclusion and
$$F_H:H_*(M)\to HF_*(O_M,\phi^1_H(O_M))
$$ is an isomorphism between singular and Floer homology groups. The construction for spectral invariants in case of conormal bundle boundary condition is done in~\cite{O1}, and in~\cite{O2} for cohomology classes. It turned out that Oh's invariants and the Viterbo's ones, are in fact the same, see~\cite{M1,M2}.

Oh proved in~\cite{O1} that these invariants are independent both on the choice of almost
complex structure $J$ (which is used in the definition of Floer homology) and, after a certain normalization,
on the choice of $H$ as far as $\phi^1_H(O_M)=L$. Using these invariants $\sigma(\alpha,L):=\sigma(\alpha,H)$, Oh derived the
non--degeneracy of Hofer's metric for Lagrangian submanifolds, the result earlier proved by Chekanov~\cite{C} using different methods. Another application to Hofer geometry is given in~\cite{M3,M4} in the characterization of geodesics in Hofer's metric for Lagrangian submanifolds of the cotangent bundle via quasi--autonomous Hamiltonians.

Spectral invariants in cotangent bundles were also studied by Monzner, Vichery and Zapolsky in~\cite{MVZ}.

\subsection{Beyond cotangent bundles}

Spectral invariants in general symplectic manifolds have been studied by several authors, and are still the subject of active research. Without attempting to give a complete references, we mention just a few. The construction of spectral invariants for contractible periodic orbits when $(P,\omega)$ is a symplectic manifold with $\omega|_{\pi_2(P)}=0$ and  $c_1|_{\pi_2(P)}=0$ was carried out by Schwarz (see~\cite{Sc}). In~\cite{L}, Leclercq constructed spectral invariants for Lagrangian Floer theory in case when $L$ is a closed submanifold of a compact (or convex in infinity) symplectic manifold $P$ and
$\omega|_{\pi_2(P,L)}=0,\quad \mu|_{\pi_2(P,L)}=0$, where $\mu$ is Maslov index. Symplectic invariants were further investigated by Eliashberg and Polterovich~\cite{EP}, Polterovich and Rosen~\cite{PR}, Oh~\cite{O4}, Humili\`ere, Leclercq and Seyfaddini ~\cite{HLS}, by Monzner, Vichery and Zapolsky~\cite{MVZ}, Lanzat~\cite{La} and also in~\cite{D},~\cite{M1,M2,M3,M4}.

\subsection{Overview of the paper} 

The above mentioned (and other) previous results concerning spectral invariants dealt either with Hamiltonian $H$ on symplectic manifold or with Lagrangian submanifold, thus they have a global character. Our result generalizes earlier constructions to the case of arbitrary open subsets of a base of cotangent bundle. We define spectral invariants in this case, and study how they intertwine with certain direct limits used in a construction.

\begin{figure}
\centering
\includegraphics[width=9cm,height=7.5cm]{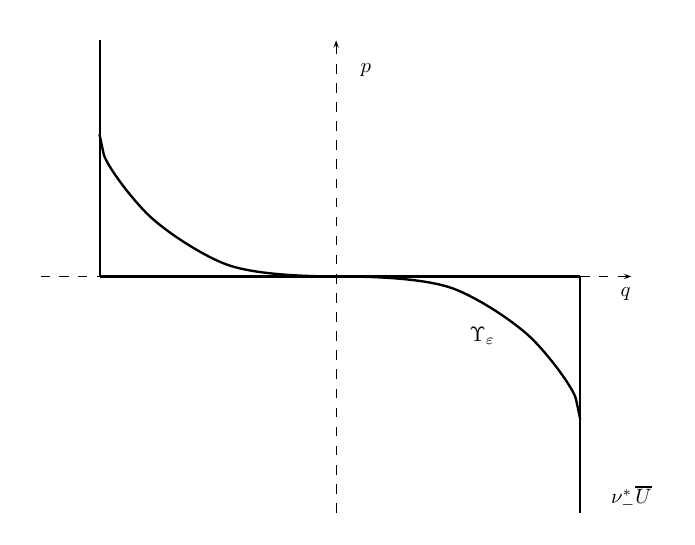}
\centering
\caption{Singular Lagrangian $\nu^*_-\overline{U}$ and approximation $\Upsilon_\varepsilon$}
\label{approx_pic}
\end{figure}

Lagrangian Floer homology for open subsets in cotangent bundles was introduced by Kasturirangan and Oh in~\cite{KO} as a part of a project of ``quantization of Eilenberg--Steenrod axioms" (see~\cite{KO1}). The construction goes as follows. Let $U\subset M$ be an open subset of a compact smooth manifold $M$, with a smooth compact
boundary $\partial U$. The conormal bundle, $\nu^*(\partial U)$, defined as
$$\nu^*(\partial U)=\{(q,p)\in T^*M\mid q\in \partial U, p|_{T_q\partial U}=0\},$$ is a Lagrangian submanifold of the cotangent bundle $T^*M$. Define
$$\nu^*_-(\partial U):=\{\alpha\in\nu^*(\partial U)\mid \alpha(\mathbf{n})\le 0,\,\mbox{for}\;\mathbf{n}\,\mbox{outward normal to}\,\partial U\}$$ and
$$\nu_-^*\overline{U}:=O_U\cup\nu^*_-(\partial U).$$ The set $\nu^*_-\overline{U}$, called the {\it (negative) conormal to $\overline{U}$}, is a singular Lagrangian submanifold, but it allows a smooth approximation by exact Lagrangian submanifolds. Let us outline a construction of these approximations, denoted by $\Upsilon_\varepsilon$, following~\cite{KO}. For $U=(-1,1)$, $M=\mathbb{R}$, $\nu^*_-\overline{U}$ and $\Upsilon_\varepsilon$ are sketched in Figure~\ref{approx_pic}.

\begin{figure}
\centering
\includegraphics[width=8cm,height=7cm]{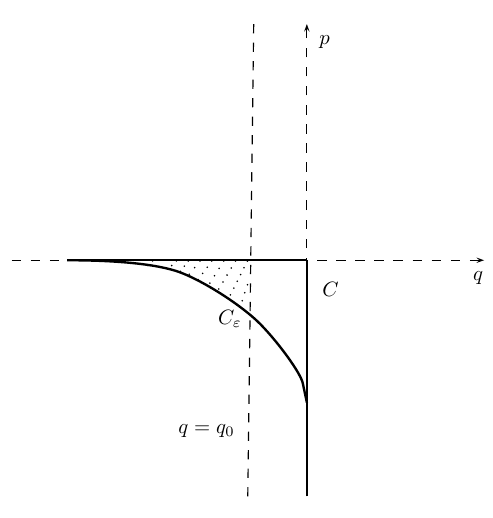}
\centering
\caption{Function $h_{\Upsilon_\varepsilon}$ is the shaded area}
\label{approx2_pic}
\end{figure}

In general case, denote by $Tb(\partial U)$ a tubular neighbourhood of $\partial U$. Since it holds:
$$Tb(\partial U)\cong\partial U\times(-1,1),$$
we have:
$$T^*M|_{Tb(\partial U)}\cong T^*(\partial U)\times((-1,1)\times\mathbb{R}).$$ Now if $C$ is a singular curve in $(-1,1)\times\mathbb{R}$:
$$C=\{(q,0)\mid -1\le q\le 0\}\cup\{(0,p)\mid p\le 0\},$$ then
$$\nu^*_-\overline{U}\cap\pi^{-1}(Tb(\partial U))=T^*(\partial U)\times C.$$ As in~\cite{KO}, denote by $C_\varepsilon$ a smooth approximation of $C$ as shown in Figure~\ref{approx2_pic} and define: 
$$\Upsilon_\varepsilon:=\nu^*_-\overline{U}\setminus\pi^{-1}(Tb(\partial U))\cup(T^*(\partial U\times C_\varepsilon).$$
To show that $\Upsilon_\varepsilon$ is exact, define a function $h_{\Upsilon_\varepsilon}:\Upsilon_\varepsilon\to\mathbb{R}$ as follows: 
\begin{itemize}
\item on $\nu^*_-\overline{U}\setminus\pi^{-1}(Tb(\partial U))=O_M|_{U\setminus Tb(\partial U)}$:  $h_{\Upsilon_\varepsilon}$ is equal to zero;
\item on the intermediate region of $\nu^*_-\overline{U}\cap\pi^{-1}(Tb(\partial U))$: $h_{\Upsilon_\varepsilon}(q_0,p_0)$ is the area of the shaded region in Figure 2 (bounded by $C_\varepsilon$, $q$-axis and the line $q=q_0$);
\item on $\nu_-^*(\partial U)\cap\Upsilon_\varepsilon$: $h_{\Upsilon_\varepsilon}$ equals to the area bounded by the $q$-axix, $p$-axes and the curve $C_\varepsilon$.
\end{itemize}

It is easy to check that $\theta|_{T\Upsilon_\varepsilon}=dh_{\Upsilon_\varepsilon}$, where $\theta$ is a canonical Liouville form on $T^*M$ and that $\Upsilon_\varepsilon\to\nu_-^*\overline{U}$ as $\varepsilon\to 0$ in Lipschitz topology (see also~\cite{KO} for more details).

Floer homology for the open set $U$ is defined to be a direct limit of Floer homologies of approximations. In order to have the latter well defined, one needs to choose a compactly supported Hamiltonian $H:T^*M\times [0,1]\to\mathbb R$ such that
 $$\phi^1_H(O_M)\pitchfork O_M$$ and
 \begin{equation}\label{prop:trans_cond}\phi^1_H(O_M)\cap O_M|_{\partial U}=\emptyset,\quad \phi^1_H(O_M)\pitchfork \nu_-^*\overline{U}.
\end{equation}
Both of the above conditions can be obtained by generic choice of $H$. Floer homology for the pair $(O_M,\Upsilon_\varepsilon)$ is now defined in a standard way, the set of the generators $CF(O_M,\Upsilon_\varepsilon:H)$ consists of the Hamiltonian paths
\begin{equation}\label{ham_paths}\dot{x}=X_H(x),\quad x(0)\in O_M,\,x(1)\in \Upsilon_\varepsilon,\end{equation} which are critical points of the {\it effective action functional:}
\begin{equation}\label{eq:ef_act_f-al}
\mathcal{A}_H^{\Upsilon_\varepsilon}(\gamma):=\int \gamma^*\theta-\int_0^1H(\gamma(t),t)dt-h_{\Upsilon_\varepsilon}(\gamma(1)).
\end{equation}

The boundary map $\partial_{J,H}$ is defined by a number of perturbed holomorphic discs with boundary on $O_M$ and $\Upsilon_\varepsilon$:
\begin{equation}\label{d}\left\{\begin{array}{l}
u:\mathbb{R}\times[0,1]\to T^*M\\
\frac{\partial u}{\partial s}+J_\varepsilon\left(\frac{\partial u}{\partial t}-X_H(u)\right)=0\\
u(s,0)\in O_M,u(s,1)\in\Upsilon_\varepsilon.
\end{array}\right.
\end{equation}
Here $J_\varepsilon$ is an almost complex structure compatible to the standard symplectic form $\omega=-d\theta$, which coincides with the canonical almost complex structure $J_0$ on $T^*M$ at infinity. By the canonical almost complex structure $J_0$ we assume the one induced by the Levi-Civita connection for a fixed metric $g_0$.

Denote by $HF_*(O_M,\Upsilon_\varepsilon:H,J_\varepsilon)$ the corresponding Floer homology (grading is given by Maslov index, see Subsection~\ref{subsec:iso_appr}).

Floer homology of the open subset $U$ is defined as a direct limit of above Floer homologies for the approximations $\Upsilon_\varepsilon$:
\begin{equation}\label{eq:dir_limit}
HF_*^-(H,U:M):=\limarr HF_*(O_M,\Upsilon_\varepsilon:H,J_\varepsilon),
\end{equation}
after defining an appropriate partial ordering to the set of pairs $(\Upsilon_\varepsilon,J_\varepsilon)$  (see Section~\ref{sec:PSS} below or~\cite{KO} for more details). The symbol $-$ in $HF_*^-$ indicates that we are dealing with the negative conormal. Defined in this way, Floer homology is isomorphic to singular homology $HF_*(U)$. More precisely,  for a special choice of Morse function $f$, such that $\nabla f$ points outward $\partial U$, Floer homology $HF^-_*(H,U:M)$ is isomorphic to Morse homology $HM(f,U)$, which, in turn, is isomorphic to $H_*(U)$.

Unlike in the paper~\cite{KO}, we have also to deal with the {\it positive conormal} to $\overline{U}$, defined as:
$$\nu_+^*\overline{U}:=O_U\cup\nu^*_+(\partial U),$$
where
$$\nu^*_+(\partial U):=\{\alpha\in\nu^*(\partial U)\mid \alpha(\mathbf{n})\ge 0,\,\mbox{for}\;\mathbf{n}\,\mbox{outward normal to}\,\partial U\}.$$
We define Floer homology $HF_*^+(H,U:M)$ in this case in the same way, as a direct limit of Floer homologies for approximations, but now, this limit will be isomorphic to the relative homology $H_*(U,\partial U)$. Again, this isomorphism is realized via Morse homology $HM_*(f,U)$, with the different choice of $f$ (with the gradient field now pointing inward at $\partial U$). The two Floer homologies $HF_*^-(H,U:M)$ and $HF_*^+(H,U:M)$ are related via Poincar\'e duality:
$$HF_*^-(H,U:M)\cong HF_{n-*}^+(\overline H,U:M)$$
(see Subsection~\ref{subsec:iso_appr} for the details).

The main aim of the paper is the construction of PSS isomorphism and the investigation of spectral invariants for the Floer homology of the open subset.

The first step in this direction is the construction of Piunikhin-Salamon-Schwarz isomorphism between $HF_*^-(H,U:M)$ (respectively $HF_*^+(H,U:M)$) and singular homology of $U$ (respectively relative homology $HF_*(U,\partial U)$) modelled by Morse homology. We will first construct PSS homomorphism for approximations. Morse homology for open subset is well defined for a fixed Morse function $f$  and a generic choice of Riemannian metric $g$, without any direct limit construction. However, in order to obtain all transversality conditions for moduli spaces of mixed type that figure in PSS homomorphisms for approximations, we have to choose (a priori) different Riemmanian metric for different $\Upsilon_\epsilon$. Therefore we will also consider Morse homology as a direct limit:
$$HM_*(f,U):=\limarr HM_*(f,U;g_s)$$
(see Section~\ref{sec:PSS}.)

More precisely, in Section~\ref{sec:PSS} we prove the following theorem.

\begin{thrm}\label{thm:PSS} Let $f^\pm{\in}\mathcal F^{\pm}(M)$ be two Morse functions from a special class of Morse functions (see Definition~\ref{defn:f} in Section~\ref{sec:PSS}). There exist PSS-type isomorphisms
$$\Phi:HM_*(f^-,U)\to HF^-_*(H,U:M),\quad
\Psi:HF^+_*(H,U:M)\to HM_*(f^+,U)$$
which are natural with respect to canonical isomorphisms in Morse and Floer theory. More precisely, if
$$\mathbf{S}_{\alpha\beta}: HF^-_*(H_\alpha,U:M)\to  HF^-_*(H_\beta,U:M),\quad  \mathbf{T}_{\alpha\beta}: HM_*(f^-_\alpha,U)\to HM_*(f^-_\beta,U)$$ denote the canonical isomorphisms in Floer and Morse theory respectively, and $\Phi_\alpha$ and $\Phi_\beta$ the corresponding PSS homomorphisms, then the diagram
$$\begin{array}{ccc}
 HF^-_*(H_\alpha,U:M) & \stackrel{\mathbf{S}_{\alpha\beta}}{\longrightarrow} &
HF^-_*(H_{\beta},U:M) \\
\Phi_\alpha\uparrow & & \uparrow\Phi_\beta \\
 HM_*(f^-_{\alpha},U) & \stackrel{\mathbf{T}_{\alpha\beta}}{\longrightarrow}
& HM_*(f^-_{\beta},U) \end{array}$$
commutes, and the same holds for the isomorphism $\Psi$.
\end{thrm}
We construct PSS homomorphisms and prove Theorem~\ref{thm:PSS} in Section~\ref{sec:PSS}. First we construct the corresponding homomorphisms for approximations $HF_*(O_M,\Upsilon_\varepsilon:H,J_\varepsilon)$ and prove
that they commute with the homomorphisms that define the direct limit~(\ref{eq:dir_limit}).

Next, we construct three pair-of-pants type products in Morse and Floer theory for open sets. Products in Morse and Floer theory were studied by various authors: Abbondandolo and Schwarz~\cite{AS}, Auroux~\cite{Aur}, Oh~\cite{O2} and also in~\cite{KMS}.

Here we establish the following products for open subset.
\begin{thrm}\label{thm:prod_open}
There exist a pair-of-pants type products:
$$\begin{aligned}
&\circ:HF_*(H_1,U:M)\otimes HF_*(H_2,U:M)\to HF_*(H_3,U:M)\\
&\cdot:HM_*(f_1,U)\otimes HM_*(f_2,U)\to HM_*(f_3,U)\\
&\star:HM_*(f,U)\otimes HF_*(H,U:M)\to HF_*(H,U:M)
\end{aligned}$$
that turns Floer homology for an open set $HF_*(H,U:M)$ into a $HM_*(f,U)-$module. The above products satisfy:
$$\Phi(\alpha\cdot\beta)=\Phi(\alpha)\circ\Phi(\beta),$$ where $\Phi$ is a PSS isomorphism from Theorem~\ref{thm:PSS}.
\end{thrm}

Theorem~\ref{thm:prod_open} is proven in Section~\ref{sec:prod}. Since Floer homology for the open set is defined as a direct limit, the key step is to prove that the products defined on homology for approximation commute with the homomorpshisms that define the direct limit.

Finally, using the above PSS isomorphism, we construct the spectral invariants for Lagrangian Floer homology of the open subset $HF(H,U:M)$.

We prove the following properties of these spectral invariants: their continuity with respects to $H$ and their subadditivity with respect to the products from Theorem~\ref{thm:prod_open}.  We also compare the above spectral invariants with the invariants for periodic orbits case, using the homomorphisms defined via ``chimneys" introduced by Abbondadolo and Schwarz~\cite{AS1}, and Albers~\cite{A}. Further, we prove the inequality of spectral invariants between two open sets
$U\stackrel{\imath}{\hookrightarrow}V$ and a specific singular homology class (see Subsection~\ref{subsec:inv_open}). This slightly generalizes a result by Oh~\cite{O3} for a spectral invariant
$$c_+(H,U):=\inf\{\lambda\in\mathbb R\mid \imath^\lambda_*:HF^\lambda_*(H,U:M)\to HF_*(H,U:M)\;\mbox{is surjective}\}.$$

More precisely, in Section~\ref{sec:inv} we prove the following theorem.

\begin{thrm}\label{thm:inv_open}
For given singular or Morse homology class $\alpha\in HM_*(f,U)\setminus\{0\}$, the spectral invariant $c_U(\alpha,H)$ defined via PSS isomorphism from Theorem~\ref{thm:PSS} has the following properties:
\begin{itemize}
\item[(A)] {\bf triangle inequality.} For $\alpha\cdot\beta\neq 0$ it holds:
$$c_U(\alpha\cdot\beta,H_1\sharp H_2)\le c_U(\alpha,H_1)+c_U(\beta,H_2)$$
\item[(B)] {\bf continuity.} relative spectral invariant $C_U(\alpha,H):=c_U(\alpha,H)-c_U(1,H)$ is continuous with respect to the Hofer norm of $H$
\item[(C)] {\bf comparison with periodic orbit invariants.} Let $\rho(\cdot,H)$ stands for a spectral invariants for periodic orbit case in $T^*M$,  $\imath_*$ is a homomorphism in homology induced by the inclusion map, and $\imath_!$ is the map obtained by inclusion map and Poincar\'e duality map:
$$\imath_!:=\PD^{-1}\circ\; \imath_*\circ\PD.$$ Suppose that the Hamiltonian $H$ satisfies the conditions from Frauenfelder-Schlenk's paper~\cite{FS} (see also Subsection~\ref{subsec:per} on page~\pageref{subsec:per}).  
Then it holds:
$$\rho(\alpha,H)\ge c_U(\imath_!(\alpha),H),\quad c_U(\alpha,H)\ge\rho(\imath_*(\alpha),H)$$
\item[(D)] {\bf invariants for subsets.} Let $U\stackrel{\imath}{\hookrightarrow}V$ be two open subset of $M$ and let
$$\jmath_{*UV}:HM_*(f,U)\to HM_*(f,V)$$ (the homomorphism induced by inclusion $\jmath:U\hookrightarrow V$) be surjective. For $\alpha\in HM_*(f,U)\setminus\{0\}$ it holds:
$$c_V(\jmath_{*UV}(\alpha),H)\le c_U(\alpha,H).$$

\item[(E)] If $H$ and $K$ are two compactly supported Hamiltonians generating the same time-one-map, i.e. $\phi_H^1=\phi_K^1$, then the corresponding invariants are the same:
$$c_U(\alpha,H)=c_U(\alpha,K),$$ so we can define $c_U(\alpha,\phi)$ for a Hamiltonian diffeomorphism $\phi$. 

\end{itemize}
\end{thrm}

\section{PSS isomorphism}\label{sec:PSS}

PSS type isomorphism was originally constructed by Piunikhin, Salamon and Schwarz~\cite{PSS} for periodic orbit case, and later adapted in~\cite{KM,A} for Lagrangian case. 

One of the nice consequences of the existence of PSS isomorphism is, for example. the commutativity of the diagram:
$$\xymatrix{
HM(f^{\alpha}) \ar[d]_{\PSS}\ar[r]&HM(f^{\beta})
  \ar[d]_{\PSS}\\
 HF(H^{\alpha})\ar[r]&HF(H^{\beta})}.$$
 In order to establish the similar naturality for several homomorphisms and operators in our case, we have to carefully investigate the subtleties related to the passing to direct limit. We start with the approximations and then pass to the limit.

\subsection{Isomorphism for approximations}\label{subsec:iso_appr}

We first establish the PSS homomorphism for approximations for negative conormal case. It follows from~(\ref{prop:trans_cond}) that all solutions of Hamiltonian equation $\dot{x}=X_H(x)$ with $x(0)\in O_M$ satisfy $x(1)\notin O_M|_{\partial U}$, so by choosing $\Upsilon$ to coincide with $\nu_-^*\overline{U}$ outside the small neighbourhood of $O_M|_{\partial U}$, we may assume that all solutions of~(\ref{ham_paths}) satisfy
\begin{equation}\label{eq:Ham_end}
x(0), x(1)\in O_M\quad \mbox{or} \quad x(0)\in O_M, \;x(1)\in \nu^*_-\overline{U}.
\end{equation}
 The grading for $x\in CF(O_M,\Upsilon:H)$ is defined to be
$$\begin{aligned}&\mu(x):=\mu_M(x)+\frac{1}{2}\dim M,\quad\mbox{for}\;x(1)\in O_U\\&\mu(x):=\mu_{\partial U}(x)+\frac{1}{2}\dim (\partial U),\quad\mbox{for}\;x(1)\in \nu^*(\partial U),
\end{aligned}$$
where $\mu_S$ is a canonically assigned Maslov index, defined for any smooth closed submanifold $S\subset M$ (see Definition 5.9 in~\cite{O1}). The dimension of the space $\mathcal M(x,y,O_M,\Upsilon:H,J)$ of perturbed holomorphic discs that satisfy~(\ref{d}) and the infinity boundary conditions:
$$u(-\infty,t)=x(t),\quad u(+\infty,t)=y(t)$$ is
$$\dim \mathcal M(x,y,O_M,\Upsilon:H,J)=\mu(y)-\mu(x)$$ for all $x,y\in CF(O_M,\Upsilon:H)$ (see~\cite{KO}).

We will consider a special class of Morse functions, as in~\cite{AD} or~\cite{La}. 
\begin{defn}\label{defn:f}
For a given Riemannian metric $g$ on $M$, let $\mathcal{F}^-(g)\subset C^{\infty}(M)$ be the set of all Morse functions $f$ on $M$ such that 
\begin{itemize}
\item $\Crit(f)\cap \overline{V}=\emptyset$, where $V$ is some neighbourhood of $\partial U$;
 \item the gradient vector field $\nabla_gf$ of $f$ is everywhere transversal to $\partial U$ and points outward $U$ along $\partial U$.\label{F^+} 
\end{itemize}
Define also 
$$\mathcal{F}^+(g):=\{f\in C^{\infty}(M)\mid -f\in\mathcal{F}^-(g)\}.$$ 
\end{defn}

\begin{figure}
\centering
\includegraphics[width=8cm,height=4cm]{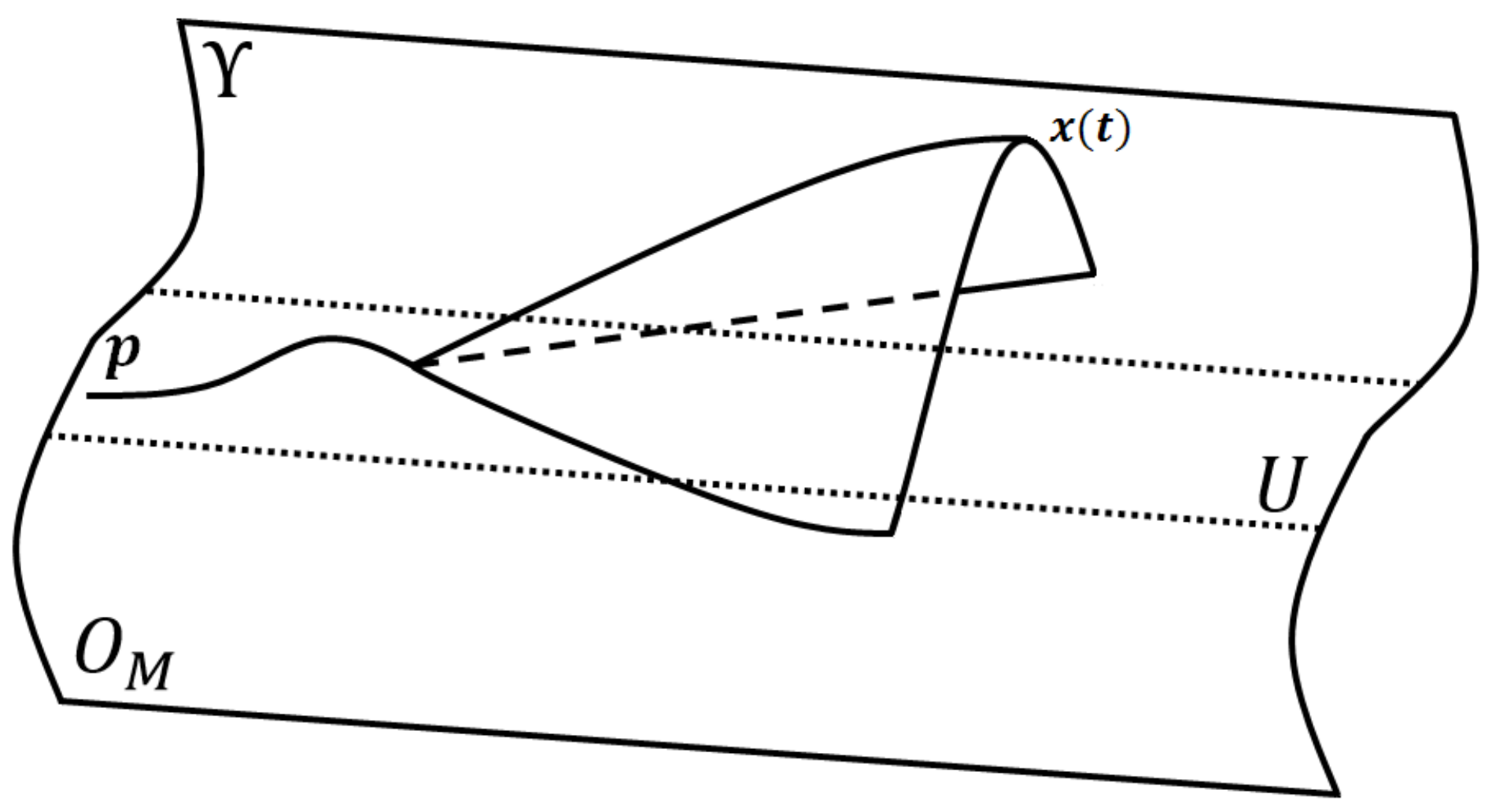}
\centering
\caption{Mixed object $\mathcal{M}(p,x)$ that defines PSS homomorphism}
\label{PSS_pic}
\end{figure}

Now let $f\in\mathcal{F}^-(g)$, $p\in\Crit(f)\cap U$ and $x\in CF(O_M,\Upsilon:H)$. Define the space of mixed objects (see Figure~\ref{PSS_pic}):
$$\begin{aligned}&\mathcal{M}(p,x):=
\mathcal{M}(p,x,O_M,\Upsilon:f,H,J,g):=\\
&\left\{(\gamma,u)\left|\begin{array}{l}
\gamma:(-\infty,0]\to U,\; u:[0,+\infty)\times[0,1]\to T^*M\\
\dot{\gamma}(s)=-\nabla_g f(\gamma(s)) \\
\frac{\partial u}{\partial s}+
J(\frac{\partial u}{\partial t}-X_{\rho_RH}(u))= 0\\
u(s,0), u(0,t)\in O_M,\,u(s,1)\in\Upsilon  \\
\gamma(-\infty)=p,\,
u(+\infty,t)=x(t) \\
u(0,1)=\gamma(0)\\
\end{array}\right.\right\}\end{aligned}$$
where $\rho_R:[0,+\infty)\rightarrow{\mathbb R}$ is a smooth function
such that
\begin{equation}\label{eq:rho}
\rho_R(s)=\begin{cases} 1, & s\ge R \\ 0, & s\le R-1. \end{cases}
\end{equation}

\medskip
\noindent Let $m_f(p)$ denotes the Morse index of a critical point $p$.

\begin{prop}\label{prop:dim} For generic choices the set ${\mathcal M}(p,x)$ is a smooth manifold
of dimension $m_f(p)-\mu(x)$.
\end{prop}
\noindent{\it Proof:} Let $x\in CF(O_M,\Upsilon:H)$ and $D$ be a half-strip $[0,1]\times[0,+\infty)$. Denote by $W^{1,r}_u(D)$ be a completion of a tangent space $T_{u}C^{\infty}(D)$ of:
\begin{equation}\label{eq:c^infty}
C^{\infty}(D):=\{u\in C^{\infty}(D,T^*M)\mid u(s,0), u(0,t)\in O_M,u(s,1)\in\Upsilon, u(+\infty,t)=x(t)\}
\end{equation}
which is
\begin{equation}\label{eq:T_u}
T_{u}C^{\infty}(D)=\left\{\eta\in C^{\infty}(D,TT^*M)\left|
\begin{array}{l}
\eta(s,t)\in T_{u(s,t)}T^*M\\
\eta([0,\infty)\times\{0\})\subset TO_M\\
\eta(\{0\}\times[0,1])\subset TO_M\\
\eta([0,\infty)\times\{1\})\in T\Upsilon\\
\eta(t,+\infty)=0
\end{array}\right.\right\}
\end{equation}
in Sobolev norm:
$$\|\eta\|_{W^{1,r}}=\left(\iint\limits_D\left(\left|\eta\right|^{r}+\left|\nabla_s\eta\right|^{r}+
\left|\nabla_t\eta\right|^{r}\right)dsdt\right)^{\frac{1}{r}}.$$
By $C^{\infty}(D,T^*M)$ in~(\ref{eq:c^infty}) and~(\ref{eq:T_u}) we mean smooth on interior of $D$ and continuous on $D$. Banach space $W^{1,r}_u(D)$ gives rise to Banach manifolds of mappings  $\mathcal{P}^{1,r}(D)$ by
$$T_{u}\mathcal{P}^{1,r}(D)=W^{1,r}_u(D).$$ 

We choose a metric $g$ as in~\cite{KO} in such a way that it becomes a product metric on a tubular neighbourhood of $\partial U$:
$$Tb(\partial U)\cong \partial U\times (-1,1).$$ Since $TT^*M$ symplectically splits into
$$TT^*M|_{Tb(\partial U)}=T(T^*(\partial U))\oplus T(T^*(-1,1)),$$ the choice of $g$ gives rise to the splitting of vertical and horizontal spaces:
$$V(TT^*M)=V((T^*(\partial U))\oplus V(T^*(-1,1)),\quad 
H(TT^*M)=H((T^*(\partial U))\oplus H(T^*(-1,1)).$$

Now let $x\in CF(O_M,\Upsilon:H)$ and $u:D\to T^*M$ the solution of
\begin{equation}\label{eq:u}
\left\{\begin{array}{l}
\frac{\partial u}{\partial s}+
J(\frac{\partial u}{\partial t}-X_{\rho_RH}(u))= 0\\
u(s,0), u(0,t)\in O_M,\,u(s,1)\in\Upsilon.
\end{array}\right.
\end{equation}
Let $n=\dim M$. We fix a trivialization 
$$\Phi_{+}:x^*(TT^*M)\to [0,1]\times \mathbb{C}^n,$$
and extend it to a trivialization
$$\Phi:u^*(TT^*M)\to D\times \mathbb{C}^n,$$ that preserves the splitting
$$\Phi(H)=\mathbb{R}^n,\quad \Phi(V)=i\mathbb{R}^n.$$ The choice of $g$ provides the splitting
$$\begin{aligned}&\Phi(V(T^*(\partial U))=i\mathbb{R}^{n-1}\times\{0\},\quad
&\Phi(V(T^*(-1,1))=\{0\}\times i\mathbb{R},\\
&\Phi(H(T^*(\partial U))=\mathbb{R}^{n-1}\times\{0\},\quad
&\Phi(H(T^*(-1,1))=\{0\}\times\mathbb{R}.\end{aligned}$$

The operator $\partial_{J,\rho H}$ defined as 
$$\partial_{J,\rho H}u=\frac{\partial u}{\partial s}+
J\left(\frac{\partial u}{\partial t}-X_{\rho_RH}(u)\right)$$ is a section of a suitable vector bundle over $\mathcal{P}^{1,r}(D)$. Denote its covariant li\-ne\-a\-ri\-za\-tion at $u$ by $L_u$. The operator $(\Phi_u)^*L_u$ is a Cauchy-Riemann type operator $\frac{\partial}{\partial s}+J\frac{\partial}{\partial t}+T$ acting on
$$W^{1,r}_\Phi:=\{\eta\in W^{1,r}(D,\mathbb{C}^n)\mid \eta(s,0),\eta(0,t)\in\mathbb{R}^n,\eta(s,1)\in\Lambda^{\Phi}(s)\},$$
where
$$\Lambda^{\Phi}(s):=\Phi(T_{u(s,1)}\Upsilon).$$

Now we proceed as in Appendix in~\cite{O2} to conclude that, for generic choice of $J$, the set $W^s(x,H)$ of $u$ satisfying~(\ref{eq:u}) is a smooth manifold of dimension $-\mu(x)+n$.

Denote by $W^u(p,f)$ the unstable manifold of
$p$. For a generic choice of parameters the evaluation map
$$\ev: W^u(p,f)\times W^s(x,H)\rightarrow U\times U, \quad
(\gamma, u)\mapsto \left(\gamma(0), u(0,1)\right)$$
is transversal to the diagonal, so
$${\mathcal M}(p,x)=\ev^{-1}(\Delta)$$ is a smooth manifold
of codimension $n$ in $W^u(p,f)\times W^s(x,H)$. Since $\dim W^u(p,f)=m_f(p)$ (see~\cite{M}) and
$\dim W^s(x,H)=-\mu(x)+n$, the dimension of ${\mathcal M}(p,x)$ is
$$\dim {\mathcal M}(p,x)=m_f(p)-\mu(x)+n-n=m_f(p)-\mu(x).$$
\qed
\bigskip

Define $\widehat{\mathcal M}(p,q)$ to be the set of all solutions of the differential equation
$$\left\{\begin{array}{l}
\gamma:\mathbb{R}\to U\\
\dot{\gamma}(s)=-\nabla_g f(\gamma(s))\\
\gamma(-\infty)=p,\;\gamma(+\infty)=q
\end{array}\right.$$ modulo $\mathbb{R}$ action and, similarly, denote by
\begin{equation}\label{eq:unpar}
\widehat{\mathcal M}(x,y):={\mathcal{M}}(x,y,O_M,\Upsilon:H,J)/\mathbb{R}.
\end{equation}

\begin{prop}\label{prop:comp} For generic choices of parameters the following is true.
\begin{itemize}
\item[(1)] If $m_f(p)=\mu(x)$, then the zero-dimensional manifold ${\mathcal M}(p,x)$ is compact, and hence, finite sets.
\item[(2)] For $m_f(p)=\mu(x)+1$, the topological boundary of the one-dimensional manifold ${\mathcal M}(p,x)$ is
$$\partial{\mathcal M}(p,x)=\bigcup_q
\widehat{\mathcal{M}}(p,q)\times\mathcal {M}(q,x)\cup
\bigcup_y
\mathcal {M}(p,y)\times\widehat{\mathcal{M}}(y,x)
$$ where the first union is taken over all $q\in\Crit(f)$,
with $m_f(q)=m_f(p)-1$, and the second over all $y\in CF(O_M,\Upsilon:H)$, such that $\mu(y)=\mu(x)+1$. 

\end{itemize}
\end{prop}

\noindent{\it Proof.} The proof follows from standard arguments, using the Arzela-Ascoli and Gromov compactness theorems.
Bubbling cannot occur due to exactness of $\omega$ and exact Lagrangian boundary conditions. Our choice of a Morse function $f\in\mathcal{F}^-(g)$ guarantees that there are no additional boundary components coming from the sequences $\gamma_n(0)=u_n(0,1)$, since it is isolated from the boundary $\partial U$.\qed

The part (1) in the previous proposition enables us to define the homomorphism between Morse and Floer homology.
Denote by
$$\begin{aligned}
CM_k(f,U):=\mathbb{Z}_2\left\langle\, p\in\Crit(f)\cap U\mid m_f(p)=k\right\rangle\\
CF_k(O_M,\Upsilon:H):=\mathbb{Z}_2\left\langle x\in CF(O_M,\Upsilon:H)\mid \mu(x)=k\right\rangle\end{aligned}$$
(i.e. $\mathbb{Z}_2-$vector spaces over the sets of generators of corresponding indices).
Let $HM_k(f,U:g)$ and $HF_k(O_M,\Upsilon:H,J)$ denote the corresponding Morse and Floer homology groups.

Denote:
$$n(p,x):=\,\sharp\,\mathcal {M}(p,x)\pmod 2$$
and define
$$\phi^{\Upsilon}:CM_k(f,U)\to CF_k(O_M,\Upsilon:H),\quad \phi^{\Upsilon}: p\mapsto\sum_{x\in CF_k(O_M,\Upsilon:H)}n(p,x)x.$$

\medskip

Before we define the homomorphisms $\psi^\Upsilon:CF_k(O_M,\Upsilon:H)\to CM_k(f,U)$, we need to describe Floer homology construction in positive conormal case.


As in~\cite{KO1}, we consider the anti-symplectic involution
\begin{equation}\label{eq:zeta}
\zeta:x=(q,p)\mapsto \overline{x}:=(q,-p).
\end{equation}
Note that $\zeta$ maps the negative conormal $\nu_-^*\overline{U}$ to the positive conormal $\nu_+^*\overline{U}$. If $\Upsilon$ is an exact Lagrangian approximation of 
$\nu_-^*\overline{U}$, then $\overline{\Upsilon}:=\zeta(\Upsilon)$ is an exact Lagrangian approximation of $\nu_+^*\overline{U}$. Next, if we define
$$\overline{H}(x,t):=-H(\zeta(x),t),\quad \overline{J}:=\zeta^*J,$$ we have 
$$CF_k(O_M,\Upsilon:H)\cong CF_{n-k}(O_M,\overline\Upsilon:\overline H).$$ We also have an identification of the space of perturbed holomorphic discs defining the boundary operation:
$$\zeta:\mathcal{M}(x,y,O_M,\Upsilon:H,J)\stackrel{\cong}{\longrightarrow}\mathcal{M}(\overline x,\overline y,O_M,\overline \Upsilon:\overline H,\overline J),$$ so $\zeta$ induces a Poincar\'e dual isomorphism:
$$\PD_F=\zeta_*:HF_k(O_M,\Upsilon:H,J)\stackrel{\cong}{\longrightarrow}
HF_{n-k}(O_M,\overline \Upsilon:\overline H,\overline J).$$

\begin{rem}
Anti-symplectic involution $\zeta$ also induces the Poincar\'e dual in Morse case, since
$$\zeta_*f=-f.$$
\end{rem}

Now choose a Morse function $f\in\mathcal{F}^+(g)$ (see Definition~\ref{defn:f}). For $\overline{x}\in CF_k(O_M,\overline{\Upsilon}:\overline{H})$, define
\begin{equation}\label{eq:mfd_comb}
\begin{aligned}&\mathcal{M}(\overline{x},p):=
\mathcal{M}(\overline{x},p,O_M,\overline{\Upsilon}:f,\overline{H},\overline{J},g):=\\
&\left\{(u,\gamma)\left|\begin{array}{l}
u:(-\infty,0]\times[0,1]\to T^*M,\;\gamma:[0,+\infty)\to U\\
\frac{\partial u}{\partial s}+
\overline{J}(\frac{\partial u}{\partial t}-X_{\tilde{\rho}_R\overline{H}}(u))= 0\\
\dot{\gamma}(s)=-\nabla_g f(\gamma(s)) \\
u(s,0), u(0,t)\in O_M,\,u(s,1)\in\overline{\Upsilon}  \\
u(-\infty,t)=\overline{x}(t),\,\gamma(+\infty)=p \\
\gamma(0)=u(0,1),\\
\end{array}\right.\right\}\end{aligned}
\end{equation}
where $\tilde{\rho}_R(s):={\rho}_R(-s)$.

As in Proposition~\ref{prop:dim} we conclude that the set $\mathcal{M}(\overline x,p)$
is a smooth manifold of dimension $\mu(\overline x)-m_f(p)$, compact in the dimension zero and with the similar description of a boundary in the dimension one:
$$\partial{\mathcal M}(\overline{x},p)=
\bigcup_{\overline{y}}
\widehat{\mathcal{M}}(\overline{x},\overline{y})\times\mathcal {M}(\overline{y},p)\cup
\bigcup_q
\mathcal {M}(\overline{x},q)\times\widehat{\mathcal{M}}(q,p).
$$
 
For $\mu(\overline x)=m_f(p)$, denote by $n(\overline{x},p):=\,\sharp\,\mathcal {M}(\overline{x},p)\pmod 2$ and define:
$$\psi^{\overline \Upsilon}:CF_k(O_M,\overline \Upsilon:\overline H)\to CM_k(f,U),\quad \psi^{\overline \Upsilon}: \overline{x}\mapsto\sum_{p\in CM_k(f,U)}n(\overline{x},p)p.$$

The proof of the following theorem follows from the standard cobordism arguments, the part (2) of the Proposition~\ref{prop:comp} and the description of $\partial\mathcal{M}(\overline{x},p)$ from above.

\begin{prop} The homomorphism $\phi^\Upsilon$ and $\psi^{\overline \Upsilon}$ induce homomorphisms
\begin{equation}\label{eq:hom_H}
\Phi^{\Upsilon}:HM_k(f^-,U:g)\to HF_k(O_M,\Upsilon:H,J)\\
\end{equation}
and
\begin{equation}\label{hom_psi}
\Psi^{\overline \Upsilon}:HF_k(O_M,\overline{\Upsilon}:\overline{H},\overline{J})\to HM_k(f^+,U:g)
\end{equation}
on the homology level, for $f^{\pm}\in\mathcal{F}^{\pm}(g)$.
\end{prop}

\medskip

If $f\in\mathcal{F}^-(g)$, then $-f\in\mathcal{F}^+(g)$, so for such $f$ we have well defined both
$$\begin{aligned}&\Phi^{\Upsilon}:HM_k(f,U:g)\to HF_k(O_M,\Upsilon:H,J)\quad\mbox{and}\\
&\Psi^{\overline \Upsilon}:HF_k(O_M,\overline{\Upsilon}:\overline{H},\overline{J})\to HM_k(-f,U:g).\end{aligned}$$

By Poincar\'e duality in Morse homology we mean the isomorphism:
\begin{equation}\label{eq:varphi}
\PD_M:HM_k(f,U:g)\stackrel{\cong}{\longrightarrow}HM_{n-k}(-f,U:g),\quad p\mapsto p.
\end{equation}

\begin{thrm}\label{thm:iso_approx}
The diagram
$$\xymatrix{
 HM_k(f,U:g) \ar[d]_{\PD_M}^\cong\ar[r]^-{\Phi^{\Upsilon}} &HF_k(O_M,\Upsilon:H,J)
  \ar[d]_{\PD_{F}}^\cong\\
 HM_{n-k}(-f,U) &HF_{n-k}(O_M,\overline \Upsilon:\overline{H},\overline{J})\ar[l]_-{\Psi^{\overline \Upsilon}}}$$
commutes and therefore, the homomorphisms $\Phi^{\Upsilon}$ and $\Psi^{\overline \Upsilon}$ are isomorphisms.
\end{thrm}

\noindent{\it Proof:} 
We need to prove 
$$\Psi^{\overline \Upsilon}\circ\PD_F\circ\Phi^\Upsilon=\PD_M.$$
For $p\in\Crit_k(f)$ it holds:
$$\Psi^{\overline \Upsilon}\circ\PD_F\circ\Phi^\Upsilon(p)
=\sum_{m_{-f}(q)=n-k}
\left(\sum_{\mu_H(x)=k}n(p,x)n(\overline{x},q)\right)q.$$
Obviously $n(\overline{x},q)=n(x,q)$, where
$$n(x,q):=\,\sharp\,\mathcal{M}(x,p,O_M,{\Upsilon}:f,H,J,g)\,\pmod 2.$$ The number $\sum\limits_{\mu_H(x)=k}n(p,x)n(x,q)$ is a cardinality of zero-dimensional manifold:
\begin{equation}\label{eq:one_comp}
\bigcup_x\mathcal{M}(p,x,O_M,\Upsilon:f,H,J,g)\times
\mathcal{M}(x,q,O_M,\Upsilon:-f,H,J,g).
\end{equation}
The rest of the proof relies on standard cobordism arguments. The manifold~(\ref{eq:one_comp}) is one component of the boundary of an auxiliary one-dimensional manifold:
$$\overline{\mathcal M}(p,q,O_M,\Upsilon;f,H,J):=
\left\{ (\gamma_{-},\gamma_{+},u,R) \left|
\begin{array}{ll}
R\in[R_0,+\infty)\\
\gamma_{-}:(-\infty,0]\rightarrow U\\
\gamma_{+}:[0,+\infty)\rightarrow U \\
u:{\mathbb R}\times [0,1]\rightarrow T^*M \\
\frac{d\gamma_{\pm}}{dt}=-\nabla f(\gamma_{\pm})\\
\frac{\partial u}{\partial s}+
J(\frac{\partial u}{\partial t}-X_{\rho_RH}(u))= 0 \\
\gamma_-(-\infty)=p, \; \gamma_+(+\infty)=q \\
u(s,0)\in O_M, \; u(s,1)\in\Upsilon\\
u(\pm\infty,t)=\gamma_{\pm}(0)
\end{array}
\right.
\right\},
$$
where $\rho_R:{\mathbb R}\to [0,1]$ is a symmetric cut-off function:
$$
\rho_R(t)=\begin{cases} 1, & |t|\le R-1 \\ 0, & |t|\ge R. \end{cases}
$$
The second boundary component is
$$\overline{\mathcal M}_{R_0}(p,q,O_M,\Upsilon;f,H,J):=
\left\{ (\gamma_{-},\gamma_{+},u) \left|
\begin{array}{ll}
\gamma_{-}:(-\infty,0]\rightarrow U\\
\gamma_{+}:[0,+\infty)\rightarrow U \\
u:{\mathbb R}\times [0,1]\rightarrow T^*M \\
\frac{d\gamma_{\pm}}{dt}=-\nabla f(\gamma_{\pm})\\
\frac{\partial u}{\partial s}+
J(\frac{\partial u}{\partial t}-X_{\rho_{R_0}H}(u))= 0 \\
\gamma_-(-\infty)=p, \; \gamma_+(+\infty)=q \\
u(s,0)\in O_M, \; u(s,1)\in\Upsilon\\
u(\pm\infty,t)=\gamma_{\pm}(0)
\end{array}
\right.
\right\}, 
$$ and the remaining components are such that induce zero mappings in homology level. This means that the mapping $\Psi^{\overline \Upsilon}\circ\PD_F\circ\Phi^\Upsilon$ is equal to
$$p\mapsto\sum_qn_{R_0}(p,q)q,$$ where $n_{R_0}(p,q)$ is a cardinality of $\overline{\mathcal M}_{R_0}(p,q,O_M,\Upsilon;f,H,J)$. Now, by standard cobordism arguments one shows that the latter mapping does not depend on $R_0$ on the homology level. Therefore we can choose $R_0=0$ and obtain holomorphic map $u$ with the boundary  on $O_M\cup \Upsilon$, so it must be constant due to the exactness of both Lagrangian submanifolds and the fact that $h_\Upsilon|_{O_U}=0$. Hence $\Psi^{\overline \Upsilon}\circ\PD_F\circ\Phi^\Upsilon$ is chain homotopic to the map obtained by counting the pairs $(\gamma_1,\gamma_2)$ with properties:
$$\left\{\begin{array}{l}
\gamma_1:(-\infty,0]\to U,\;\gamma_2:[0,+\infty)\to U\\
\dot{\gamma}_j=-\nabla f(\gamma_j)\\
\gamma_1(-\infty)-p,\;\gamma_2(+\infty)=q\\
\gamma_1(0)=\gamma_2(0).
\end{array} \right.$$ The trajectory $\gamma_1\sharp\gamma_2$ is a negative gradient trajectory of $f$ connecting two critical points of the same Morse index. Number of such pairs is equal 1 in case $p=q$ and 0 otherwise. Therefore, $\Psi^{\overline \Upsilon}\circ\PD_F\circ\Phi^\Upsilon$ is chain homotopic to the homomorphism $\PD_M$. 
\qed
\bigskip

For two Morse functions $f_\alpha,f_\beta\in\mathcal{F}^{\pm}$, Morse homologies
$HM(f_\alpha,U:g)$ and $HM(f_\beta,U:g)$ are canonically isomorphic (see~\cite{Sc1}). Similarly, for two Hamiltonians $H_\alpha$ and $H_\beta$, the corresponding Floer homologies $HF(O_M,\Upsilon:H_\alpha,J)$ and $HF(O_M,\Upsilon:H_\beta,J)$ are isomorphic (see~\cite{KO}). Denote these canonical isomorphisms by
$$\begin{aligned}
T_{\alpha\beta}:HM(f_\alpha,U:g)&\stackrel{\cong}{\longrightarrow}HM(f_\beta,U:g)\\
S^\Upsilon_{\alpha\beta}:HF(O_M,\Upsilon:H_\alpha,J)&\stackrel{\cong}{\longrightarrow}HF(O_M,\Upsilon:H_\beta,J).
\end{aligned}
$$ Note that we use the same notation, $S_{\alpha\beta}$ and $T_{\alpha\beta}$, for canonical isomorphisms for two Morse homologies (relative and absolute one, i.e. for Morse functions from both $\mathcal{F}^+(g)$ and $\mathcal{F}^-(g)$) and two Floer homologies (negative and positive conormal case). 

Denote by
$$\begin{array}{l}
\Phi^\Upsilon_\alpha:HM_k(f_{\alpha},U:g)\to HF_k(O_M,\Upsilon:H_{\alpha},J)\\
\Psi^\Upsilon_\alpha:HF_k(O_M,\Upsilon:H_{\alpha},J)\to HM_k(f_{\alpha},U:g)\end{array}
$$
the homomorphisms defined in~(\ref{eq:hom_H}) and~(\ref{hom_psi}).

\begin{thrm}\label{thm:comm_approx} The diagrams
\begin{equation}\label{diag:S:T}
\begin{array}{ccc}
 HF_k(O_M,\overline{\Upsilon}:\overline H_\alpha,\overline J) & \stackrel{S^{\overline \Upsilon}_{\alpha\beta}}{\longrightarrow} &
HF_k(O_M,\overline \Upsilon:\overline H_{\beta},\overline J) \\
\Psi^{\overline \Upsilon}_\alpha\downarrow & & \downarrow\Psi^{\overline \Upsilon}_\beta \\
 HM_k(f_{\alpha},U:g) & \stackrel{T_{\alpha\beta}}{\longrightarrow}
& HM_k(f_{\beta},U:g)
\end{array}
\end{equation} 
and
$$\begin{array}{ccc}
 HF_k(O_M,\Upsilon:H_\alpha,J) & \stackrel{S^\Upsilon_{\alpha\beta}}{\longrightarrow} &
HF_k(O_M,\Upsilon:H_{\beta},J) \\
\Phi^\Upsilon_\alpha\uparrow & & \uparrow\Phi^\Upsilon_\beta \\
 HM_k(f_{\alpha},U:g) & \stackrel{T_{\alpha\beta}}{\longrightarrow}
& HM_k(f_{\beta},U:g)
\end{array}$$
commute.
\end{thrm}

\noindent{\it Proof:} The homomorphism $T_{\alpha\beta}\circ\Psi^{\overline \Upsilon}_{\alpha}$ is the same as the map $K$ defined on generators as
$$K(\overline x_\alpha):=\sum_{p_\beta}\tilde n(\overline x_{\alpha},p_{\beta}),$$ where
$\tilde n(\overline x_{\alpha},p_{\beta})$ is the cardinal number (modulo $2$) of zero-dimensional component of the smooth manifold

$$\widetilde{\mathcal M}_T(\overline x_{\alpha},p_{\beta},O_M,\overline \Upsilon:\overline H_{\alpha},f_{\alpha\beta,T},\overline J):=
\left\{ (\gamma,u) \left|
\begin{array}{l}
u:(-\infty,0]\times[0,1]\rightarrow T^*M \\
\gamma:[0,+\infty)\rightarrow U \\
\frac{\partial u}{\partial s}+
\overline J(\frac{\partial u}{\partial t}-X_{\overline\rho_R\overline H_{\alpha}}(u))= 0 \\
u(-\infty,t)=\overline x_\alpha(t)\\
u(s,0), u(0,t)\in O_M,\,u(s,1)\in\overline \Upsilon \\
\dot\gamma(s)=-\nabla f_{\alpha\beta,T}(\gamma) \\
\gamma(+\infty)=p_\beta\\
\gamma(0)=u(0,1).
\end{array}
\right.
\right\}
$$
Here $f_{\alpha\beta,T}(s)\in\mathcal{F}^+(g)$ satisfies
$$f_{\alpha\beta,T}(s)=\begin{cases} f_{\alpha}, & s\le T \\ f_{\beta}, & s\ge 2T.\end{cases}$$
for fixed $T>0$, and $\overline\rho_R(s)=\rho_R(-s)$, for $\rho_R$ defined in~(\ref{eq:rho}).

Indeed, to see this, consider the boundary of one-dimensional auxiliary manifold
$$\widetilde{\mathcal M}(\overline x_{\alpha},p_{\beta},O_M,\overline \Upsilon:\overline H_{\alpha},f_{\alpha\beta,T},\overline J):=
\left\{ (\gamma,u,T) \left|
\begin{array}{l}
u:(-\infty,0]\times[0,1]\rightarrow T^*M \\
\gamma:[0,+\infty)\rightarrow U \\
\frac{\partial u}{\partial s}+
\overline J(\frac{\partial u}{\partial t}-X_{\overline\rho_R\overline H_{\alpha}}(u))= 0 \\
u(-\infty,t)=\overline x_\alpha(t)\\
u(s,0), u(0,t)\in O_M,\,u(s,1)\in\overline \Upsilon \\
\dot\gamma(s)=-\nabla f_{\alpha\beta,T}(\gamma) \\
\gamma(+\infty)=p_\beta\\
\gamma(0)=u(0,1).
\end{array}
\right.
\right\}
$$

Similarly, $\Psi^{\overline \Upsilon}_{\beta}\circ S^{\overline \Upsilon}_{\alpha\beta}$ is the same as the map
$$L(\overline x_\alpha):=\sum_{p_\beta}\check{n} (\overline x_{\alpha},p_{\beta}),$$ where
$\check n(\overline x_{\alpha},p_{\beta})$ is the number of zero-dimensional component of the smooth manifold

$$\check{\mathcal M}_T(\overline x_{\alpha},p_{\beta},O_M,\overline \Upsilon:H_{\alpha\beta,T},f_{\alpha},\overline J):=
\left\{ (\gamma,u) \left|
\begin{array}{l}
u:(-\infty,0]\times[0,1]\rightarrow T^*M \\
\gamma:[0,+\infty)\rightarrow U \\
\frac{\partial u}{\partial s}+
\overline J(\frac{\partial u}{\partial t}-X_{\overline\rho_R H_{\alpha\beta,T}}(u))= 0 \\
u(-\infty,t)=\overline x_\alpha(t)\\
u(s,0), u(0,t)\in O_M,\,u(s,1)\in\overline \Upsilon \\
\dot\gamma(s)=-\nabla f_{\alpha}(\gamma) \\
\gamma(+\infty)=p_\beta\\
\gamma(0)=u(0,1).
\end{array}
\right.
\right\}
$$
where $H_{\alpha\beta,T}(s)$ is Hamiltonian function satisfying
$$H_{\alpha\beta,T}(s)=\begin{cases} \overline H_{\alpha}, & s\le -2T \\ \overline H_{\beta}, & s\le -T.\end{cases}.$$

So we need to proof that the maps $K$ and $L$ are the same in homology level.

Fix $T>0$. Let 
$(f_{\alpha\beta}^{\lambda},H_{\alpha\beta}^{\lambda})_{0\le\lambda\le 1}$
be a homotopy connecting
$(f_{\alpha\beta}^{\lambda},H_{\alpha\beta}^{\lambda})|_{\lambda=0}=
(f_{\alpha},H^s_{\alpha\beta,T})$ and
$(f_{\alpha\beta}^{\lambda},H_{\alpha\beta}^{\lambda})|_{\lambda=1}=
(f_{\alpha\beta,T}^s,\overline H_{\beta})$

Let $p_\beta\in U$ be a critical point of $f_{\beta}$ and $\overline x_{\alpha}\in CF(O_M,\overline \Upsilon:\overline H_{\alpha})$. 
Define the auxiliary $(m_{f_\beta}(p_\beta)-\mu_{\overline{H}_\alpha}(\overline x_{\alpha})+1)-$dimensional manifold:
$$
\widehat{\mathcal M}(\overline x_{\alpha},p_{\beta},O_M,\overline \Upsilon:H_{\alpha\beta}^{\lambda},f_{\alpha\beta}^{\lambda},\overline J):=
\left\{ (\gamma,u,\lambda) \left|
\begin{array}{l}
u:(-\infty,0]\times[0,1]\rightarrow T^*M \\
\gamma:[0,+\infty)\rightarrow U \\
\frac{\partial u}{\partial s}+
\overline J(\frac{\partial u}{\partial t}-X_{\rho_RH_{\alpha\beta}^{\lambda}}(u))= 0 \\
u(-\infty,t)=\overline x_\alpha(t)\\
u(s,0), u(0,t)\in O_M,\,u(s,1)\in\overline \Upsilon \\
\dot\gamma(s)=-\nabla f_{\alpha\beta}^{\lambda}(\gamma) \\
\gamma(+\infty)=p_\beta\\
\gamma(0)=u(0,1).
\end{array}
\right.
\right\}
$$
For $m_{f_\beta}(p_\beta)=\mu_{\overline{H}_\alpha}(\overline x_{\alpha})$ the boundary of one dimensional manifold $\widehat{\mathcal M}(\overline x_{\alpha},p_{\beta},O_M,\overline \Upsilon:H_{\alpha\beta}^{\lambda},f_{\alpha\beta}^{\lambda},\overline J)$ is
$$\begin{aligned}&\bigcup_{\overline y_{\alpha}}\widehat{\mathcal{M}}(\overline x_{\alpha},\overline y_{\alpha})\times\widehat{\mathcal M}(\overline y_{\alpha},p_{\beta},O_M,\overline \Upsilon:H_{\alpha\beta}^{\lambda},f_{\alpha\beta}^{\lambda},\overline J)\;\cup\\
&\bigcup_{q_{\beta}}\widehat{\mathcal M}(\overline x_{\alpha},q_{\beta},O_M,\overline \Upsilon:H_{\alpha\beta}^{\lambda},f_{\alpha\beta}^{\lambda},\overline J)\times
\widehat{\mathcal{M}}(q_{\beta},p_{\beta})\;\cup\\
&\mathcal{M}(\overline x_{\alpha},p_{\beta},O_M,\overline \Upsilon:f_{\alpha\beta,T},\overline H_{\alpha},\overline J,g)\;\cup\;
\mathcal{M}(\overline x_{\alpha},p_{\beta},O_M,\overline \Upsilon:f_{\beta},H_{\alpha\beta,T},\overline J,g).
\end{aligned}
$$ 
The rest of the proof relies on standard arguments, see e.g.~\cite{PSS,KM}.
\qed

\subsection{Isomorphism for Floer homology of open set}

In order to define Floer homology for the open set as a direct limit of Floer homologies for the approximations, Kasturirangan and Oh defined a partial ordering on the set of approximations as:
$$\Upsilon_a\le\Upsilon_b\;\Longleftrightarrow\;\varphi_a\le \varphi_b\;\mbox{on}\;U.$$
The function $\varphi_a$ is defined by $h_a=\varphi_a\circ\pi$ on $U$, where $h_a:\Upsilon_a\to \mathbb{R}$ is a smooth function such that $\theta|_{T\Upsilon_a}=dh_a$ (recall that $\Upsilon_a$ is exact) and $\pi:T^*M\to M$ is a canonical projection. Since $H$ is fixed, one has to vary the almost complex structure $J$ to obtain a generic condition for Fredholm theory. Denote by $J_a$ an almost complex structure corresponding to $\Upsilon_a$ and denote by
$$\mathbf F_{ab}:HF_k(O_M,\Upsilon_a:H,J_a)\to HF_k(O_M,\Upsilon_b:H,J_b)$$
a canonical homomorphism that satisfies:
$$\mathbf F_{ac}=\mathbf F_{bc}\circ \mathbf F_{ab}$$
for given triple $\Upsilon_a\le\Upsilon_b\le\Upsilon_c$ sufficiently close to $\nu^*\overline{U}$ (see~\cite{KO}). As we have mentioned in Introduction, Floer homology for an open subset $U$, modelled by negative conormal, is defined as
$$HF_k^-(H,U:M):=\limarr HF_k(O_M,\Upsilon_s:H,J_s).$$

Since we want to establish an isomorphism between Floer homology and Morse homology for a fixed Morse function, we will vary Riemannian metric, so the term ``generic choices" in the Proposition~\ref{prop:trans_cond} refers to an almost complex structure $J$ and Riemannian metric $g$.

Fix a Hamiltonian function $H$ and a Morse function $f$. For a Lagrangian approximation $\Upsilon_a$, choose an almost complex structure $J_a$ and a Riemannian metric $g_a$ such that all the transversality conditions are fulfilled, i.e. the sets $\widehat{\mathcal{M}}(p,q)$, $\widehat{\mathcal{M}}(x,y)$, $\mathcal{M}(p,x)$ and $\mathcal{M}(x,p)$ are manifolds for all $p,q\in\Crit(f)$ and all Hamiltonian paths $x,y$ with boundaries on $O_M$ and $\Upsilon_a$. For two Riemannian metric $g_a$ and $g_b$ there is a canonical isomorphism
$$\mathbf G_{ab}:HM_k(f,U:g_a)\to HM_k(f,U:g_b)$$
satisfying
$$\mathbf G_{ac}=\mathbf G_{bc}\circ \mathbf G_{ab},\quad G_{aa}=\Id.$$ This functoriality allows to consider the set $\{HM_*(f,U:g_a)\}$ as a directed system and to define Morse homology $HM_k(f,U)$ as a direct limit:
$$HM_k(f,U):=\limarr HM_k(f,U:g_a):=
\bigsqcup_sHM_k(f,U:g_s)/\sim$$ 
where
$$p_a\sim p_b\Leftrightarrow\mathbf{F}_{ac}(p_a)=\mathbf{F}_{bc}(p_b)$$ for some $c$. The set $HM_k(f,U)$ obviously has a vector space structure and is isomorphic to all $HM_k(f,U:g)$.

Consider a diagram:
\begin{equation}\label{eq:diag_ups_g}
\begin{array}{lllllllll}
\cdots&\longrightarrow&HM_k(g_a)&\stackrel{\mathbf G_{ab}}{\longrightarrow}&HM_k(g_b)&\stackrel{\mathbf G_{bc}}{\longrightarrow}&HM_k(g_c)&\longrightarrow&\cdots\\
&&\downarrow\Phi^a&&\downarrow\Phi^b&&\downarrow\Phi^c&&\\
\cdots&\longrightarrow&HF_k(\Upsilon_a)&\stackrel{\mathbf F_{ab}}{\longrightarrow}&HF_k(\Upsilon_b)&\stackrel{\mathbf F_{bc}}{\longrightarrow}&HF_k(\Upsilon_c)&\longrightarrow&\cdots
\end{array}
\end{equation}
where we use the abbreviations
$$\begin{aligned}
&HM_k(g_a):=HM_k(f,U:g_a)\\
&\Phi^a:=\Phi^{\Upsilon_a}\\
&HF_k(\Upsilon_a):=HF_k(O_M,\Upsilon_a:H,J_a),
\end{aligned}
$$
and so on.

\begin{prop}\label{prop:var}
The diagram~(\ref{eq:diag_ups_g}) commutes.
\end{prop}
\noindent{\it Proof:} The homomorphism $G_{ab}$ at the chain level (we denoted by $\mathbf G_{ab}$ the induced homomorphism in homology) is defined via the cardinal number of the set
\begin{equation}\label{eq:G_ab}
\mathcal{M}(p,q:\tilde{g}_s):=
\left\{ \gamma \left|
\begin{array}{l}
\gamma:\mathbb{R}\rightarrow U \\
\frac{d\gamma}{ds}=-\nabla_{\tilde g_s} f(\gamma) \\
\gamma(-\infty)=p,\,\gamma(+\infty)=q,\end{array}\right.\right\}
\end{equation}
and the homomorphism $F_{ab}$ via the number of elements in
\begin{equation}\label{eq:F_ab}
\mathcal{M}(x,y:\widetilde{\Upsilon}_s):=
\left\{ u\left|
\begin{array}{l}
u:\mathbb{R}\times[0,1]\rightarrow T^*M\\
\frac{\partial u}{\partial s}+
\widetilde J_s(\frac{\partial u}{\partial t}-X_H(u))= 0 \\
u(s,0)\in O_M,\,u(s,1)\in\widetilde{\Upsilon}_s\\
u(-\infty,t)=x(t),\;u(+\infty,t)=y(t).
\end{array}\right.\right\}
\end{equation}
Here:
\begin{itemize}
\item $\widetilde\Upsilon_s$ is a monotone homotopy for $s\in\mathbb{R}$ such that
\begin{equation}\label{eq:mon_hom}
\widetilde\Upsilon_s=\begin{cases}
\Upsilon_a, &s\le-R\\
\Upsilon_b, &s\ge R;
\end{cases}\end{equation} 
(by monotone homotopy we mean $s_1\le s_2\Rightarrow\Upsilon_{s_1}\le\Upsilon_{s_2}$)
\item $\widetilde J_s$ is a corresponding family of generic almost complex structures;
\item $\tilde g_s$ is a homotopy of Riemannian metrics such that
$$\tilde g_s=\begin{cases}
g_a, &s\le-T\\
g_b, &s\ge T.
\end{cases}$$
\end{itemize}

The rest proof of Proposition~\ref{prop:var} relies on cobordism arguments, similarly to the proof of Theorem~\ref{thm:comm_approx}, so we omit the details.\qed

\bigskip

We have the similar partial ordering for the set of approximations of positive conormal. Actually, we define such a partial ordering via anti-symplectic involution:
$$\overline \Upsilon^a\le\overline \Upsilon^b\Leftrightarrow\zeta(\Upsilon^a)\le\zeta(\Upsilon^b),$$ where $\zeta$ is defined in~(\ref{eq:zeta}). We define Floer homology for an open subset $U$, modelled by positive conormal, as
$$HF_k^+(H,U:M):=\limarr HF_k(O_M,\overline \Upsilon_s:\overline H,\overline J_s).$$

Let $\mathbf{F}^+_{ab}$ denote the canonical isomorphism for the positive conormal:
$$\mathbf F^+_{ab}:HF_k(O_M,\overline\Upsilon_a:\overline H,\overline J_a)\to HF_k(O_M,\overline \Upsilon_b:\overline H,\overline J_b)$$
defined in the same way as $\mathbf{F}_{ab}$, by the number of solutions of~(\ref{eq:F_ab}) (see also~\cite{KO}).

\begin{thrm}
There exist direct limit homomorphisms
\begin{equation}\label{PSS_denote}
\Phi:HM_k(f^-,U)\rightarrow HF_k^-(H,U:M).
\end{equation}
and
$$
\Psi:HF^+_k(H,U:M)\rightarrow HM_k(f^+,U).$$\qed
\end{thrm}

\noindent{\it Proof.} The diagram
\begin{equation}\label{eq:sim_comm_diag}
\begin{array}{lllllllll}
\cdots&\longrightarrow&HF_k(\overline\Upsilon_a)&\stackrel{\mathbf{F}^+_{ab}}{\longrightarrow}&HF_k(\overline\Upsilon_b)&\stackrel{\mathbf{F}^+_{bc}}{\longrightarrow}&HF_k(\overline\Upsilon_c)&\longrightarrow&\cdots\\
&&\downarrow\Psi^a&&\downarrow\Psi^b&&\downarrow\Psi^c&&\\
\cdots&\longrightarrow&HM_k(g_a)&\stackrel{\mathbf G_{ab}}{\longrightarrow}&HM_k(g_b)&\stackrel{\mathbf G_{bc}}{\longrightarrow}&HM_k(g_c)&\longrightarrow&\cdots
\end{array}
\end{equation}
commutes. This can be proved in the same way as Proposition~\ref{prop:var}. Now the proof follows directly from Proposition~\ref{prop:var} and the commutative diagram~(\ref{eq:sim_comm_diag}).
\qed

\bigskip

The Poincar\'e duality isomorphism $\PD_M$ defined in~(\ref{eq:varphi}) obviously commutes with the maps $\mathbf{G}_{ab}$, being defined as $p\mapsto p$. Hence it induces an isomorphism on a direct limit Morse homology $HM(f,U)$. Denote it again by 
$$\PD_M:HM_k(f,U)\stackrel{\cong}{\longrightarrow}HM_{n-k}(-f,U).$$
In order to emphasize the particular Riemannian metric we will use the notation:
$$\PD_M^a:HM_k(f,U:g_a)\stackrel{\cong}{\longrightarrow}HM_{n-k}(-f,U:g_a).$$

Regarding the Floer case, it is easy to see that
$$\mathbf F^+_{ab}=\PD_F\circ\mathbf{F}_{ab}\circ\PD_F^{-1},
$$ 
so $\PD_F$ defines the map  
$$\PD_F:HF_k^-(H,U:M)\stackrel{\cong}{\longrightarrow}HF_{n-k}^+(\overline H,U:M).$$
Again, denote:   
$$\PD_F^a:HF^-_k(O_M,\Upsilon^a:H,J_a)\stackrel{\cong}{\longrightarrow}HF^+_{n-k}(O_M,\overline\Upsilon^a:\overline H,\overline J_a).$$

\begin{thrm}
The diagram
$$\xymatrix{
 HM_k(f,U) \ar[d]_{\PD_M}^\cong\ar[r]^-{\Phi} &HF_k^-(H,U:M)
  \ar[d]_{\PD_{F}}^\cong\\
 HM_{n-k}(-f,U) &HF_{n-k}^+(\overline{H},U:M)\ar[l]_-{\Psi}}$$
commutes and therefore, the induced maps $\Phi$ and $\Psi$ are isomorphisms. 
\end{thrm}

\noindent{\it Proof:}
From Theorem~\ref{thm:iso_approx} we have
$$
\Psi^a\circ\PD_F^a\circ\Phi^a=\PD_M^a.
$$
Let $p_a\in HM_k(f,U:g_a)$ be the representative of the class $[p_a]\in HM_k(f,U)$. We have
$$\begin{aligned}&\Psi\circ\PD_F\circ\Phi([p_a])=\Psi\circ\PD_F([\Phi^a(p_a)])=\Psi([\PD_a\circ\Phi^a(p_a))]=\\&[\Psi^a\circ\PD_F^a\circ\Phi^a(p_a)]=[\PD_M^a(p_a)]=\PD_M([p_a]).\end{aligned}$$

\qed
\bigskip

From the canonical isomorphisms
$$S^a_{\alpha\beta}:=S^{\Upsilon_a}_{\alpha\beta}:HF(O_M,\Upsilon_a:H_\alpha,J_a)\stackrel{\cong}{\longrightarrow}HF(O_M,\Upsilon_a:H_\beta,J_a)$$ and the commutativity of the diagrams
$$\begin{array}{lllllllll}
\cdots&\longrightarrow&HF_k(\Upsilon_a:H_\alpha)&\stackrel{F^\alpha_{ab}}{\longrightarrow}&HF_k(\Upsilon_b:H_\alpha)&\stackrel{F^\alpha_{bc}}{\longrightarrow}&HF_k(\Upsilon_c:H_\alpha)&\longrightarrow&\cdots\\
&&\downarrow S^a_{\alpha\beta} &&\downarrow S^b_{\alpha\beta}&&\downarrow S^c_{\alpha\beta}&&\\
\cdots&\longrightarrow&HF_k(\Upsilon_a:H_\beta)&\stackrel{F^\beta_{ab}}{\longrightarrow}&HF_k(\Upsilon_b:H_\beta)&\stackrel{F^\beta_{bc}}{\longrightarrow}&HF_k(\Upsilon_c:H_\beta)&\longrightarrow&\cdots
\end{array}$$
(and similarly for the positive conormal) we obtain isomorphisms:
$$\begin{aligned}&\mathbf{S}_{\alpha\beta}^-:HF_k^{-}(H_\alpha,U:M)\stackrel{\cong}{\longrightarrow}HF_k^-(H_{\beta},U:M)\\
&\mathbf{S}_{\alpha\beta}^+:HF_k^{+}(\overline H_\alpha,U:M)\stackrel{\cong}{\longrightarrow}HF_k^+(\overline H_{\beta},U:M).\end{aligned}$$
Similarly, we have
$$\mathbf{T}_{\alpha\beta}:HM_k(f_\alpha,U)\stackrel{\cong}{\longrightarrow}HM_k(f_{\beta},U).$$
\begin{thrm} The diagram
$$\begin{array}{ccc}
 HF_k^+(\overline H_\alpha,U:M) & \stackrel{\mathbf{S}^+_{\alpha\beta}}{\longrightarrow} &
HF_k^+(\overline H_{\beta},U:M) \\
\Psi_\alpha\downarrow & & \downarrow\Psi_\beta \\
 HM_k(f_{\alpha},U) & \stackrel{\mathbf{T}_{\alpha\beta}}{\longrightarrow}
& HM_k(f_{\beta},U)  \end{array}$$
commutes and the same holds for the other PSS isomorphisms, $\Phi_\alpha$ and $\Phi_\beta$.
\end{thrm}
\noindent{\it Proof:} Recall that the diagram~(\ref{diag:S:T}) commutes for all approximations close enough to $\nu_+^*\overline{U}$ and for generic choices. So we have
$$\begin{aligned}&\mathbf{T}_{\alpha\beta}\circ
\Psi_\alpha([x_a])=
\mathbf{T}_{\alpha\beta}\left(\left[\Psi_\alpha^a(x_a)\right]\right)=
\left[T^a_{\alpha\beta}\left(\Psi_\alpha^a(x_a)\right)\right]=\\
&[\Psi_\beta^a\left(S^a_{\alpha\beta}(x_a)\right)]=\Psi_\beta\left([S^a_{\alpha\beta}(x_a)]\right)=\Psi_\beta\circ
\mathbf{S}_{\alpha\beta}([x_a]),\end{aligned}$$
for every $[x_a]\in HF_k^+(\overline H_\alpha,U:M)$.
\qed

This proves Theorem~\ref{thm:PSS}.

\bigskip

\section{Product on homology and module structure}\label{sec:prod}

In this section we construct a product $\circ$ on Floer homology for an open subset, a product $\cdot$ on Morse homology, and a product $\star$ which turns Floer homology to a module over a Morse homology ring. We also prove the compatibility of PSS isomorphisms with the above product and thus we prove Theorem~\ref{thm:prod_open}. 

\subsection{Product on homology}\label{subsec:product}
\begin{figure}
\centering
\includegraphics[width=8cm,height=3cm]{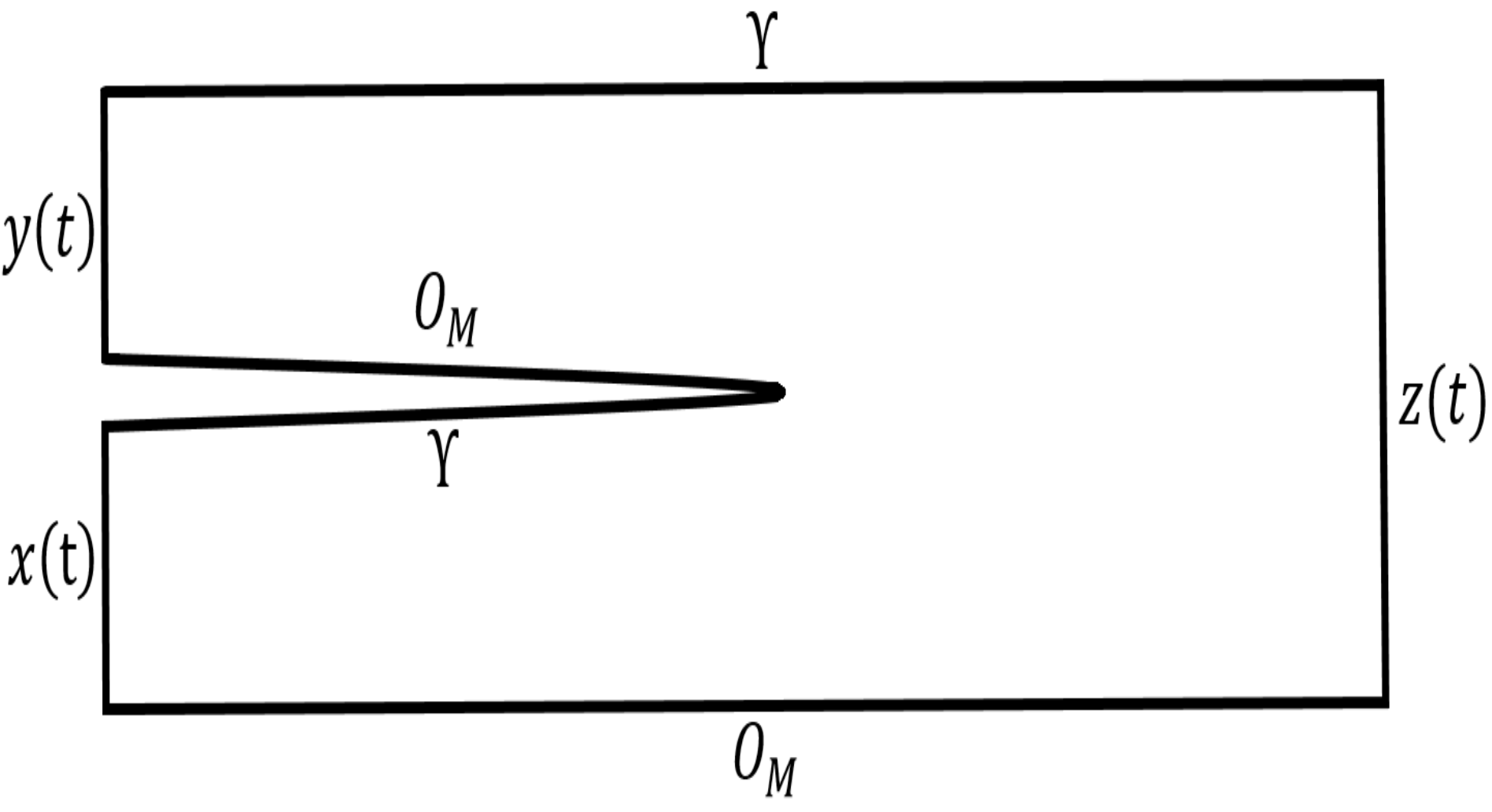}
\centering
\caption{Moduli space ${\mathcal M}(x,y;z)$ that defines a product $\circ$}
\label{Pants_pic}
\end{figure}

First we construct a product on Floer homology for an open subset
$$
\circ:HF_*(H_1,U:M)\otimes HF_*(H_2,U:M)\longrightarrow HF_*(H_3,U:M).
$$ In order to do that, we need to define a product $\circ$ on homology for approximation
\begin{equation}\label{product*}
\circ:HF_*(O_M,\Upsilon:H_1,J_\Upsilon)\otimes HF_*(O_M,\Upsilon:H_2,J_\Upsilon)\longrightarrow HF_*(O_M,\Upsilon:H_3,J_\Upsilon),
\end{equation}
and to check its compatibility with direct limit homomorphisms. A product~(\ref{product*}) is defined by a number of pair--of--pants objects. More precisely, let $\Sigma$ be a Riemannian surface (with a boundary)\label{SigmaUpsilon}
$$
{\mathbb R}\times[-1,0]\sqcup{\mathbb R}\times[0,1]
$$
with the identification $(s,0^-)\sim(s,0^+)$ for $s\geq0$ (see Figure~\ref{Pants_pic}).\\

\medskip

Denote by $\Sigma_1$, $\Sigma_2$, $\Sigma_3$ the three ends
 $$\Sigma_j\approx [0,1]\times(-\infty,0]$$
and by $u_j:=u|_{\Sigma_j}$, $j=1,2,3$.
Let $\rho:\mathbb{R}\to[0,1]$ denote the smooth cut--off function such that
$$\rho(s)=\begin{cases}1,&s\le -2\\0,&s\ge -1.\end{cases}$$
For $x\in CF_*(O_M,\Upsilon:H_1)$, $y\in CF_*(O_M,\Upsilon:H_2)$ and $z\in CF_*(O_M,\Upsilon:H_3)$ we define a moduli space
$$\begin{aligned}&{\mathcal M}(x,y;z)=&\left\{u:\Sigma\to T^*M\left|\begin{array}{l}
 \partial_su_j+J_\Upsilon(\partial_tu_j-X_{\rho_jH_j}\circ u_j)=0,\,j=1,2,3\\
 \partial_su+J_\Upsilon\partial_tu=0,\,\mbox{on}\;\Sigma_0:=\Sigma\setminus(\Sigma_1\cup \Sigma_2\cup \Sigma_3)\\
 u(s,-1)\in O_M,\,u(s,1)\in\Upsilon,\,s\in{\mathbb R}\\
 u(s,0^-)\in\Upsilon,\,u(s,0^+)\in O_M,\,s\le 0\\
 u_1(-\infty,t)=x(t)\\
 u_2(-\infty,t)=y(t)\\
 u_3(-\infty,t)=z(t)
 \end{array}\right.\right\}\end{aligned}$$
 (see Figure~\ref{Pants_pic}).

For generic choices, ${\mathcal M}(x,y;z)$ is a smooth $(\mu(x)+\mu(y)+\mu(z)-2n)$--dimensional manifold. For two generators $x$ and $y$ of Floer homology, a map $\circ$ is defined as
$$
x\circ y:=\sum_z\sharp_2{\mathcal M}(x,y;z)z,
$$
where, $\sharp_2{\mathcal M}(x,y;z)$ denotes the (modulo 2) number of elements of a zero-dimensional component of ${\mathcal M}(x,y;z)$. We extend the product $\circ$ to 
$$\circ:CF_*(O_M,\Upsilon:H_1)\otimes CF_*(O_M,\Upsilon:H_2)\to CF_*(O_M,\Upsilon:H_3)$$ by bilinearity. By standard cobordism arguments, one can show that $\circ$ commutes with boundary maps and induces a product in homology~(\ref{product*}).

\medskip

\medskip

The following lemma provides the compatibility of the product $\circ$ with the direct limit homomorphisms. Recall that we denote by $\mathbf{F}_{ab}$ the homomorphism
$$\mathbf{F}_{ab}: HF_*(O_M,\Upsilon_a:H,J_a)\to HF_*(O_M,\Upsilon_b:H,J_b)$$ defined by~(\ref{eq:F_ab}). Here $J_a=J_{\Upsilon_a}$, etc. To emphasize the Hamiltonian, we will write $\mathbf{F}_{ab}^{H}$.

\begin{lem}\label{lem}
For $x_a\in HF_*(O_M,\Upsilon_a:H_1,J_a)$, $y_a\in HF_*(O_M,\Upsilon_a:H_2,J_a)$ it holds
\begin{equation}\label{eq:F+circ}
\mathbf F^{H_3}_{ab}(x_a\circ y_a)=\mathbf F^{H_1}_{ab}(x_a)\circ \mathbf F_{ab}^{H_2}(y_a).\end{equation}
\end{lem}

\begin{figure}
\centering
\includegraphics[width=13cm,height=7cm]{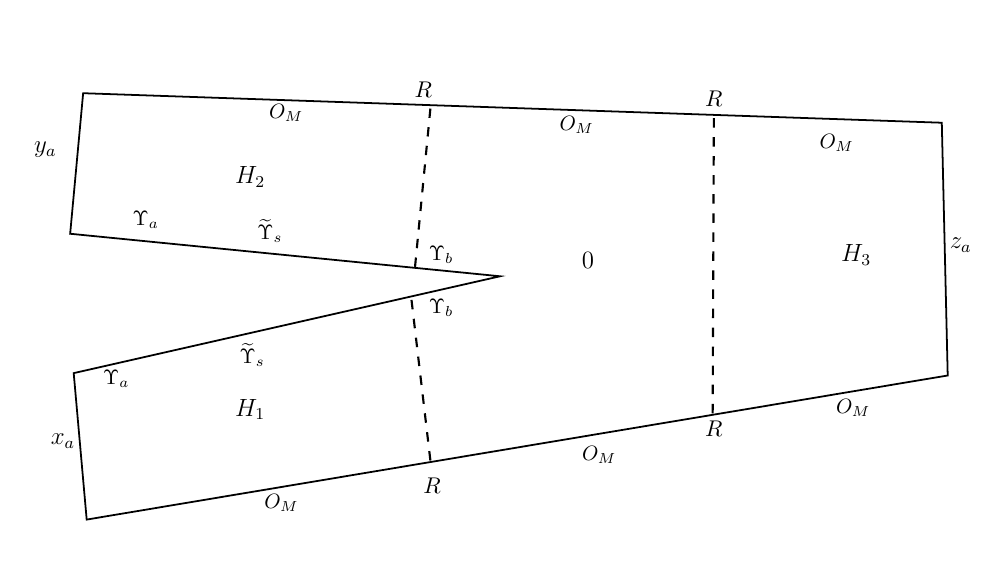}
\centering
\caption{Manifold $\mathcal{M}(x_a,y_a,z_a:\widetilde{\Upsilon}_s)$}
\label{aux_mfd_pic}
\end{figure}

\noindent{\it Proof.} The homomorphism $\mathbf{F}_{ab}$ is an isomorphism for $a$, $b$ large enough. The inverse homomorphism is actually $\mathbf{F}_{ba}$ (defined as in~(\ref{eq:F_ab}), despite the reversed order of $a$ and $b$). This can be proved using exactly the same cobordism arguments similar to ones in the proof of the independence of Floer homology with respect to the parameters (Hamiltonian, almost complex structure). Therefore,~(\ref{eq:F+circ}) is equivalent to
\begin{equation}\label{eq:F^{-1}}
x_a\circ y_a=\left(\mathbf F^{H_3}_{ab}\right)^{-1}\left(
\mathbf F^{H_1}_{ab}(x_a)\circ \mathbf F_{ab}^{H_2}(y_a)\right)=\mathbf F^{H_3}_{ba}\left(
\mathbf F^{H_1}_{ab}(x_a)\circ \mathbf F_{ab}^{H_2}(y_a)\right).
\end{equation}
In order to prove~(\ref{eq:F^{-1}}), consider the following auxiliary one-dimensional manifold. Let $\widetilde{\Upsilon}_s$ be as in~(\ref{eq:mon_hom}) and $x_a$, $y_a$, $z_a$ be the solutions of
$$\begin{aligned}
&\dot x_a(t)=X_{H_1}(x_a(t)),\quad x_a(0)\in O_M,\quad x_a(1)\in\Upsilon_a\\
&\dot y_a(t)=X_{H_2}(y_a(t)),\quad y_a(0)\in O_M,\quad y_a(1)\in\Upsilon_a\\
&\dot z_a(t)=X_{H_3}(z_a(t)),\quad z_a(0)\in O_M,\quad z_a(1)\in\Upsilon_a.
\end{aligned}$$
For $R>0$, define $\mathcal{M}_R(x_a,y_a,z_a:\widetilde{\Upsilon}_s)$ 
 to be the set of all solutions $u:\Sigma\to T^*M$ of the equation
$$\bar{\partial}_{\tilde{J},\tilde{H}}u=0$$ where $\tilde{H}$ is depicted in the Figure~\ref{aux_mfd_pic}, as well as corresponding boundary conditions. The almost complex structure $\tilde{J}$ is chosen to satisfy all the regularity conditions.  

Define $\mathcal{M}(x_a,y_a,z_a:\widetilde{\Upsilon}_s)$
to be the set of all pairs $(R,u)$, where $R\in[R_0,+\infty)$ and
$u\in\mathcal{M}_R(x_a,y_a,z_a:\widetilde{\Upsilon}_s)$.   

Now the boundary of one dimensional component of $\mathcal{M}(x_a,y_a,z_a:\widetilde{\Upsilon}_s)$ is the union of the following five strata (recall $\widehat{\mathcal{M}}(x,y)$ denotes the space of unparametrized trajectories defining the boundary operator $\partial$ in Floer homology, see~(\ref{eq:unpar}) and $\mathcal{M}(x,y;\widetilde{\Upsilon}_s)$ is defined in~(\ref{eq:F_ab})):
$$\begin{aligned}
&\mathcal{B}_1=\bigcup_{\tilde{x}_a}\widehat{\mathcal{M}}(x_a,\tilde{x}_a)
\times\mathcal{M}(\tilde{x}_a,y_a,z_a:\widetilde{\Upsilon}_s) \\
&\mathcal{B}_2=\bigcup_{\tilde{y}_a}\widehat{\mathcal{M}}(\tilde{y}_a,y_a)\times\mathcal{M}(x_a,\tilde{y}_a,z_a:\widetilde{\Upsilon}_s)\\
&\mathcal{B}_3=\bigcup_{\tilde{z}_a}\mathcal{M}(x_a,y_a,\tilde{z}_a:\widetilde{\Upsilon}_s)\times\widehat{\mathcal{M}}(\tilde{z}_a,z_a)\\
&\mathcal{B}_4=\mathcal{M}_{R_0}(x_a,y_a,z_a:\widetilde{\Upsilon}_s) \\
&\mathcal{B}_5=\bigcup_{x_b,y_b,z_b}\mathcal{M}(x_a,x_b;\widetilde{\Upsilon}_s)\times\mathcal{M}(y_a,y_b;\widetilde{\Upsilon}_s)\times\mathcal{M}_{R_1}(x_a,y_a,z_b:\widetilde{\Upsilon}_s)\times\mathcal{M}(z_b,z_a;\widetilde{\Upsilon}_s)\\
\end{aligned}$$ 
The operations induced by the number of elements of boundary strata $\mathcal{B}_1$, $\mathcal{B}_2$ and $\mathcal{B}_3$ are zero in the homology, and the operations defined by the cardinality of $\mathcal{B}_4$ and of $\mathcal{B}_5$ are equal to $x_a\circ y_a$ and $\mathbf F^{H_3}_{ba}\left(
\mathbf F^{H_1}_{ab}(x_a)\circ \mathbf F_{ab}^{H_2}(y_a)\right)$ on the homology level. \qed

\bigskip

Now we are able to define the product $\circ$ on the direct limit homology group.

\begin{prop}\label{prop:circ_open} The product $\circ$ defines a product on Floer homology for open subset:
$$
\circ:HF_*(H_1,U:M)\otimes HF_*(H_2,U:M)\longrightarrow HF_*(H_3,U:M).
$$
\end{prop}

\noindent{\it Proof.} Let $[x]\in HF_*(H_1,U:M)$ and $[y]\in HF_*(H_2,U:M)$ be the classes of elements $x\in HF_*(O_M,\Upsilon_a:H_1,J_a)$ and $y\in HF_*(O_M,\Upsilon_{a'}:H_2,J_{a'})$ in a direct limit. In general, $a$ and $a'$ are not the same, but, since $x$ and $F_{a\tilde{a}}(x)$ represent the same element in $HF_*(H_1,U:M)$ we can take $F_{a\max\{a,a'\}}(x)$ and $F_{a'\max\{a,a'\}}(y)$ as representatives of $[x]$ and $[y]$ respectively. Therefore we can assume that $x$ and $y$ belong to some $HF_*(O_M,\Upsilon:H_1,J)$ and $HF_*(O_M,\Upsilon:H_2,J)$, for the same $\Upsilon$. We now define a product $\circ$ in homology as
$$[x]\circ[y]:=[x\circ y].
$$
We need to check that a product does not depend on representatives of a class. Let $x_a$ and $x_b$ represent the same element in $HF_*(H_1,U:M)$ and similarly $y_a$ and $y_b$ in $HF_*(H_2,U:M)$. This means that there exist homomorphisms
$$\begin{aligned}
&\mathbf F_{ac}:(O_M,\Upsilon_a:H_1,J_a)\to (O_M,\Upsilon_c:H_1,J_c)\\
&\mathbf F_{bc}:(O_M,\Upsilon_b:H_1,J_a)\to (O_M,\Upsilon_c:H_1,J_c)\\
&\mathbf F_{ad}:(O_M,\Upsilon_a:H_2,J_a)\to (O_M,\Upsilon_d:H_2,J_d)\\
&\mathbf F_{bd}:(O_M,\Upsilon_b:H_2,J_a)\to (O_M,\Upsilon_d:H_2,J_d)
\end{aligned}$$
such that
$$\mathbf F_{ac}(x_a)=\mathbf F_{bc}(x_b),\quad
\mathbf F_{ad}(y_a)=\mathbf F_{bd}(y_b).$$
Let $e=\max\{c,d\}$. We have
$$\begin{aligned}
&\mathbf F_{ae}(x_a\circ y_a)\stackrel{(\ref{eq:F+circ})}{=}\mathbf F_{ae}(x_a)\circ\mathbf F_{ae}(y_a)=\mathbf F_{ce}\left(\mathbf F_{ac}(x_a)\right)\circ\mathbf F_{de}\left(\mathbf F_{ad}(y_a)\right)=\\
&\mathbf F_{ce}\left(\mathbf F_{bc}(x_b)\right)\circ\mathbf F_{de}\left(\mathbf F_{bd}(y_b)\right)=
\mathbf F_{be}(x_b)\circ \mathbf F_{be}(y_b)\stackrel{(\ref{eq:F+circ})}{=}\mathbf{F}_{be}(x_b\circ y_b)
\end{aligned}$$
which means that $x_a\circ y_a$ and $x_b\circ y_b$ represent the same element in $HF_*(H_3,U:M)$.\qed

\medskip

\subsection{Morse homology ring}

Let us recall the construction of the homology product on $HM_*(f,U:g_\Upsilon)$. Let $f_1,f_2, f_3\in\mathcal{F}^-(g)$ be three Morse functions such that $W^u_{f_1}(p_1)\pitchfork W^u_{f_2}(p_2)\pitchfork W^u_{f_3}(p_3)$ for every critical point $p_k$ of $f_k$.
For $p_i\in CM_*(f_i,U)$, $i=1,2,3$, we define the moduli space $\mathcal M(p_1,p_2;p_3)$ to be the set of all trees $\gamma:=(\gamma_1,\gamma_2,\gamma_3)$ such that
$$\left\{\begin{array}{l}\gamma_j:(-\infty,0]\to U,\,j=1,2,3\\
\dot\gamma_j=-\nabla f_j(\gamma_j),\,j=1,2,3\\
\gamma_i(-\infty)=p_i,\,i=1,2,3\\
\gamma_1(0)=\gamma_2(0)=\gamma_3(0).
\end{array}\right.$$

\begin{figure}
\centering
\includegraphics[width=6cm,height=2.5cm]{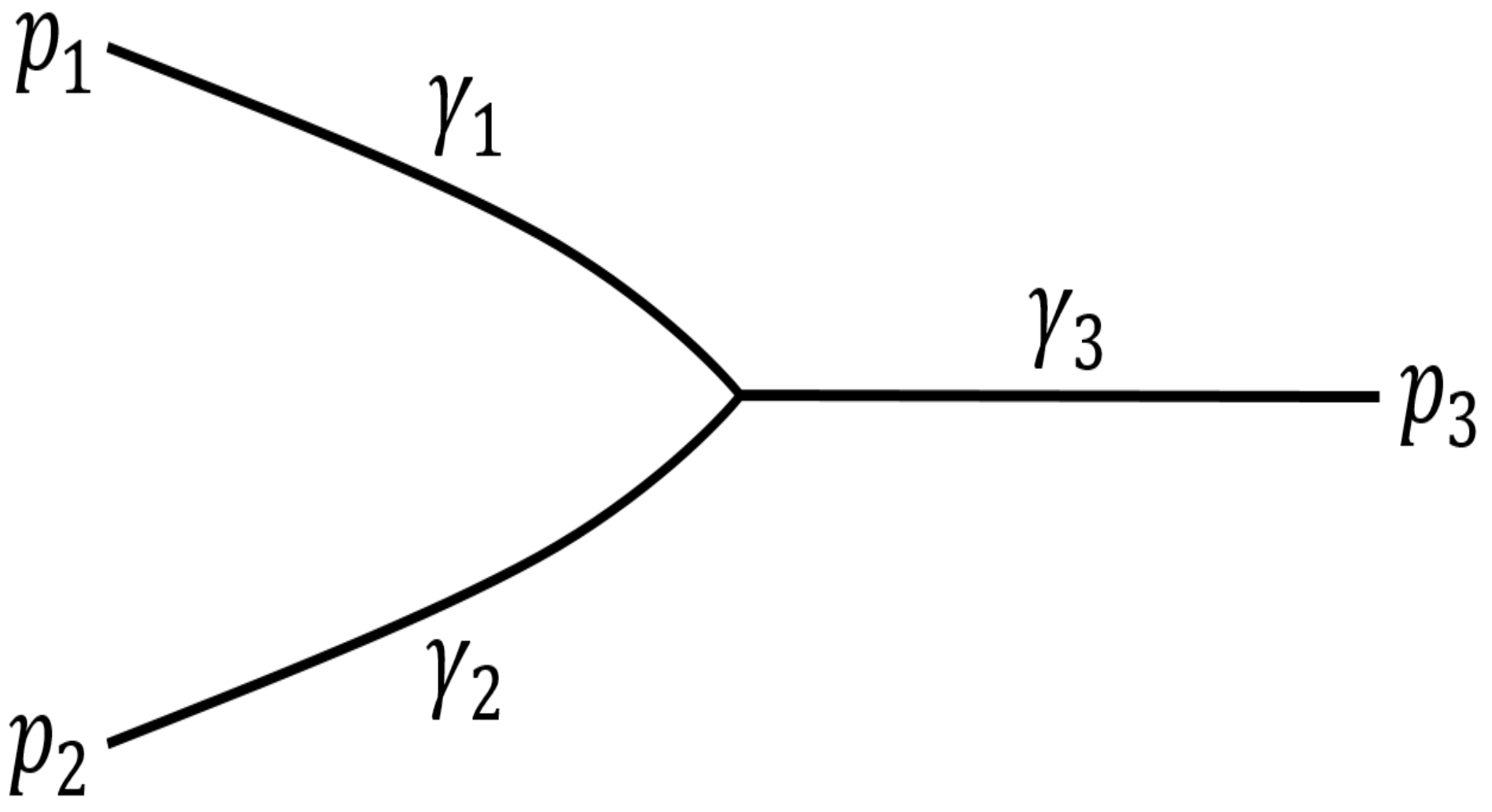}
\centering
\caption{The set of trees $\mathcal M(p,q;r)$}
\label{Tree_pic}
\end{figure}

\noindent For generic choice of $g$ these spaces are manifolds of dimension
$$
m_{f_1}(p_1)+m_{f_2}(p_2)+m_{f_3}(p_3)-2n.$$
If $n(p_1,p_2;p_3)$ denotes the mod 2 number of a zero--dimensional component, then the product $\cdot$ is defined at the chain level:
$$\cdot:CM_*(f_1,U)\otimes CM_*(f_2,U)\longrightarrow CM_*(f_3,U),
$$
as:
$$p_1\cdot p_2:=\sum_{p_3} n(p_1,p_2;p_3)p_3$$
on generators. The choice of Morse functions (see the Definition~\ref{defn:f} of $\mathcal{F}^{-}(g)$) provides that the loss of compactness of  $\mathcal M(p_1,p_2;p_3)$ is possible only as the breaking of trajectories inside $U$. Therefore $\cdot$ commutes with the Morse boundary operator and it is well defined at the homology level:
$$\cdot:HM_*(f_1,U:g)\otimes HM_*(f_2,U:g)\longrightarrow HM_*(f_3,U:g),
$$

It is also well defined as a product on a direct limit homologies:
$$\cdot:HM_*(f_1,U)\otimes HM_*(f_2,U)\longrightarrow HM_*(f_3,U),
$$
since it holds:
\begin{equation}\label{eq:morse_prod}
\mathbf{G}_{ab}^{f_3}(p_a\cdot q_a)=\mathbf{G}_{ab}^{f_1}(p_a)\cdot\mathbf{G}_{ab}^{f_2}(q_a).
\end{equation}
The latter equality can be proved in the similar way as Lemma~\ref{lem}.

The following proposition establishes the ring structure PSS isomorphism.

\begin{prop}\label{prop:circ_cdot_PSS} Let $f_1,f_2,f_3\in\mathcal{F}^-(g)$. For $[\alpha]\in HM_k(f_1,U)$, $[\beta]\in HM_k(f_2,U)$ it holds
$$\Phi_3([\alpha\cdot\beta])=\Phi_1([\alpha])\circ\Phi_2([\beta]),$$
where $\Phi_j$ is the PSS isomorphism~(\ref{PSS_denote}) obtained by Morse function $f_j$.
\end{prop}

\begin{figure}
\centering
\includegraphics[width=14cm,height=5.5cm]{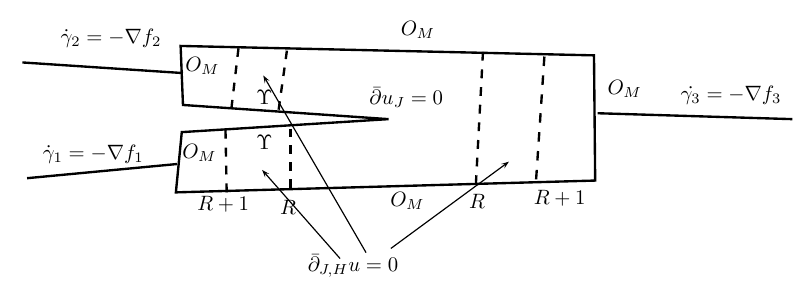}
\centering
\caption{Auxiliary one-dimensional manifold from the proof of Proposition~\ref{prop:circ_cdot_PSS}}
\label{aux_mfd2_pic}
\end{figure}

\noindent{\it Proof:}
It follows from the definition of $\Phi$ and Proposition~\ref{prop:circ_open} that it is enough to show that
\begin{equation}\label{eq:circ_cdot_app}
\Phi^{\Upsilon}_3(\alpha\cdot\beta)=\Phi_1^{\Upsilon}(\alpha)\circ\Phi_2^{\Upsilon}(\beta)\end{equation}
for a fixed approximation $\Upsilon$ and fixed Riemannian metric defining the product $\cdot$. The equality~(\ref{eq:circ_cdot_app}) is equivalent to
$$(\Phi^{\Upsilon}_3)^{-1}\left(\Phi_1^{\Upsilon}(\alpha)\circ\Phi_2^{\Upsilon}(\beta)\right)=\alpha\cdot\beta$$ and, by Theorem~\ref{thm:iso_approx} the latter equality is equivalent to
\begin{equation}\label{eq:aux:inv}
\PD_M^{-1}\circ\Psi_3^{\overline{\Upsilon}}\circ\PD_F\left(\Phi_1^{\Upsilon}(\alpha)\circ\Phi_2^{\Upsilon}(\beta)\right)=\alpha\cdot\beta
\end{equation}
The equality~(\ref{eq:aux:inv}) follows from cobordism arguments similar to ones used in the proof of Theorem~\ref{thm:iso_approx}, Proposition~\ref{prop:var} and Lemma~\ref{lem}. The auxiliary one-dimensional manifold we use here is
explained by Figure~\ref{aux_mfd2_pic}.

\qed

\bigskip

\subsection{Module structure}

Let $f\in\mathcal{F}^-(g)$. For every approximation $\Upsilon$, we can define an external product
$$
\star:CM_*(f,U:g)\otimes CF_*(O_M,\Upsilon:H,J_\Upsilon)\longrightarrow CF_*(O_M,\Upsilon:H,J_\Upsilon)
$$
by a number of a suitable mixed-type objects. More precisely, let $H^s$ denotes a smooth family of Hamiltonians such that
$$
H^s(\cdot,t)=\begin{cases}H(\cdot,1-t),&s\le -2\\0,&-1\le s\le 1\\H(\cdot,t),&s\ge2.\end{cases}
$$
For $p\in CM_*(f,U)$, $x,y\in CF_*(O_M,\Upsilon:H,J_\Upsilon)$ let ${\mathcal M}(p,x;y)$ be a moduli space of pairs $(\gamma,u)$ such that
$$\begin{aligned}&\left\{\begin{array}{l}
 \gamma:(-\infty,0]\to U,\,u:{\mathbb R}\times[0,1]\to T^*M\\
 \dot{\gamma}=-\nabla f(\gamma(t))\\
 \partial_su+J_\Upsilon(\partial_tu-X_{H^s}(u))=0\\
 u(s,0)\in \Upsilon,\,u(s,1)\in O_M,\,s\le 0\\
 u(s,0)\in O_M,\,u(s,1)\in\Upsilon,\,s\ge 0\\
 \gamma(-\infty)=p\\
 u(-\infty,t)=x(t),\,u(+\infty,t)=y(t)\\
 \gamma(0)=u(0,0).
 \end{array}\right.\end{aligned}$$
The dimension of ${\mathcal M}(p,x;y)$ equals to
$$
\mu(x)-\mu(y)+m_f(p)-n,
$$
and the zero--dimensional component is compact. Now define a product $\star$ on the set of the generators of chain complexes as:
$$
p\star x=\sum_y\sharp_2{\mathcal M}(p,x;y)y,
$$
where $\sharp_2{\mathcal M}(p,x;y)$ denotes the$\mod 2$ cardinality of the zero--dimensional component of ${\mathcal M}(p,x;y)$. Using standard cobordism arguments, as above, one can show that $\star$ induces a product in homology. Similarly to~\cite{L} one shows that 
\begin{equation}\label{eq:star_mod}
(p\cdot q)\star x=p\star(q\star x),
\end{equation}
for all $p,q\in HM_*(f,U:g)$ and $x\in HF_*(O_M,\Upsilon:H,J_\Upsilon)$.

In order to have the products $\cdot$ and $\star$ well defined on a direct limit of Morse and Floer homology groups, we need to check their compatibilities with homomorphisms $\mathbf{G}$ and $\mathbf{F}$.

\begin{lem}
Let $\mathbf{F}_{ab}$ and $\mathbf{G}_{ab}$ be the homomorphisms that define the direct limit Morse and Floer homology groups, obtained by the number of~(\ref{eq:F_ab}) and~(\ref{eq:G_ab}) respectively. Then it holds
\begin{equation}\label{eq:star_com}
\begin{aligned}
&\mathbf{G}_{ab}(p_a\cdot q_a)=\mathbf{G}_{ab}(p_a)\cdot\mathbf{G}_{ab}(q_a)\\
&\mathbf{F}_{ab}(p_a\star x_a)=\mathbf{G}_{ab}(p_a)\star\mathbf{F}_{ab}(x_a)
\end{aligned}
\end{equation}
for all $p_a,q_a\in HM_*(f,U:g_{\Upsilon_a})$ and $x_a\in HF_*(O_M,\Upsilon_a:H,J_{\Upsilon_a})$.
\end{lem}

The proof is similar to the proof of Lemma~\ref{lem}.

It follows from~(\ref{eq:star_mod}) and~(\ref{eq:star_com}) that $\cdot$ and $\star$ are well defined operations on $HM_*(f,U)$ and $HF_*^-(H,U;M)$ and that $HF_*^-(H,U;M)$ is a $HM_*(f,U)$--module.

This proves Theorem~\ref{thm:prod_open}.

\section{Spectral invariants}\label{sec:inv}

In this section we define spectral invariants for open subset and prove their properties listed in Theorem~\ref{thm:inv_open}. We define spectral invariants via PSS isomorphism constructed in Section~\ref{sec:PSS}, but they can be defined alternatively, as a limit of spectral invariants for the approximations (see Proposition~\ref{prop:limit} below). This alternative definition of spectral invariants will be the key ingredient in the proof of some properties from Theorem~\ref{thm:inv_open}.

\subsection{Invariants for the open subset}\label{subsec:inv_open}

In the rest of the paper we will only consider Floer homology for approximations and open set in the negative conormal case, as well as the corresponding PSS isomorphisms
$$\begin{aligned}&\Phi^\Upsilon:HM_k(f,U:g_\Upsilon)\stackrel{\cong}{\longrightarrow}HF_k(O_M,\Upsilon:H,J_\Upsilon)\\
&\Phi:HM_k(f,U)\stackrel{\cong}{\longrightarrow}HF_k^-(H,U:M),
\end{aligned}$$ for $f\in\mathcal{F}^-(g)$.
Therefore, we will omit the sign $-$ in $HF_k^-(H,U:M)$, in order to simplify notations.

If we consider $\mathcal{A}_H^\Upsilon$ restricted to
$$\Omega(O_M,\Upsilon):=\{\gamma\in C^{\infty}([0,1],T^*M)\mid \gamma(0)\in O_M,\,\gamma(1)\in\Upsilon\},$$
 we have
 $$d\mathcal{A}_H^\Upsilon(\gamma)(\xi)=\int_0^1\left(\omega(\dot{\gamma},\xi)-dH(\gamma)(\xi)\right)dt.$$

Recall that the {\it filtered Floer homology groups} for approximations are defined as homology groups of the filtered chain complex
$$CF^\lambda_k(O_M,\Upsilon:H):=\{x\in CF_k(O_M,\Upsilon:H)\mid \mathcal{A}^\Upsilon_H(x)<\lambda\}.$$

Since the action functional decreases along the strips that define the boundary operator
$$\partial_{J,H}:CF_k(O_M,\Upsilon:H)\to CF_{k-1}(O_M,\Upsilon:H),$$
the boundary operator descends to $CF^{\lambda}_k(O_M,\Upsilon:H)$ and defines
$$\partial_{J,H}^\lambda:CF_k^\lambda(O_M,\Upsilon:H)\to CF^\lambda_{k-1}(O_M,\Upsilon:H).$$
Denote the corresponding homology groups by $HF_k^\lambda(O_M,\Upsilon:H,J_\Upsilon)$.

Now denote by
$$\imath^{\lambda}_{\Upsilon*}:HF_k^\lambda(O_M,\Upsilon:H,J_\Upsilon)\to HF_k(O_M,\Upsilon:H,J_\Upsilon)$$ the homomorphism
induced by the inclusion map $\imath_\Upsilon^{\lambda}$ and, for $\alpha\in HM_k(f,U:g_\Upsilon)\setminus\{0\}$ define
$$c_\Upsilon(\alpha,H):=\inf\{\lambda\mid\Phi^\Upsilon(\alpha)\in\IM(\imath^\lambda_{\Upsilon*})\}.$$

We need to defined the filtered Floer homology for an open set. Recall that the direct limit homomorphisms $\mathbf{F}_{ab}$ are defined via the monotone family $\widetilde{\Upsilon}_s$ that connects $\Upsilon_a$ and $\Upsilon_b$ (see~(\ref{eq:mon_hom})). Proposition 3.4 from~\cite{KO} states that the corresponding action functional $\mathcal{A}_H^{\widetilde{\Upsilon}_s}$ decreases along perturbed holomorphic strips that define $\mathbf{F}_{ab}$, in particular, that
$$\mathcal{A}_H^{\Upsilon_b}(u(y))\le \mathcal{A}_H^{\Upsilon_a}(u(x))$$ whenever there exists an $u\in\mathcal{M}(x,y:\widetilde{\Upsilon}_s)$. Therefore the homomorphisms $\mathbf{F}_{ab}$ descend to the filtered chain complex. By standard arguments one shows that they are also well defined on filtered homology groups:
$$\mathbf{F}_{ab}^\lambda:HF_k^{\lambda}(O_M,\Upsilon_a:H,J_a)\to HF_k^{\lambda}(O_M,\Upsilon_b:H,J_b).$$
Now we define the filtered Floer homology for an open set as a direct limit:
$$HF_k^{\lambda}(H,U:M):=\limarr HF_k^{\lambda}(O_M,\Upsilon_s:H,J_s).$$ One easily verifies that
$$\mathbf{F}_{ab}^\lambda\circ\imath^\lambda_{\Upsilon_a*}=\imath^\lambda_{\Upsilon_b*}\circ\mathbf{F}_{ab}^\lambda,$$
where 
$$\imath^\lambda_{\Upsilon_a*}:HF_k^{\lambda}(O_M,\Upsilon_a:H,J_a)\to HF_k(O_M,\Upsilon_a:H,J_a)$$ denotes the inclusion-induced map for the approximations. Hence the induced inclusion maps
$$\imath^\lambda_*:HF_k^{\lambda}(H,U:M)\to HF_k(H,U:M)$$ are also well defined.

\begin{defn} Let $\alpha\in HM_k(f,U)\setminus\{0\}$. A {\it spectral invariant for an open set} is defined as
\begin{equation}\label{eq:def_inv_c_open}
c_U(\alpha,H):=\inf\{\lambda\mid \Phi(\alpha)\in\IM(\imath_*^\lambda)\}.
\end{equation}
\end{defn}

\bigskip
The natural question that occurs is the question of the relation of the spectral invariants for an open subset with the spectral invariants for the approximations, i.e. weather $c_\Upsilon(\cdot,H)$ converges to $c_U(\cdot,H)$ as $\Upsilon\to\nu^*\overline{U}$. Actually, a stronger property holds.

\begin{prop}\label{prop:limit} Let $\alpha\in HM_*(f,U:g_\Upsilon)\setminus\{0\}$. Then there exists an approximation $\widetilde{\Upsilon}$ such that
$$
c_U([\alpha],H)=c_{\overline{\Upsilon}}(G_{\Upsilon\overline{\Upsilon}}(\alpha),H)
$$
for all $\overline{\Upsilon}\le\widetilde{\Upsilon}$.
\end{prop}

\noindent{\it Proof:} We have the following commutative diagram
\begin{equation}\label{eq:2diagrams}
\begin{array}{lllllllll}
\cdots&\rightarrow&HM_k(f,U:g_{\Upsilon^a})&\stackrel{\mathbf G_{ab}}{\rightarrow}&HM_k(f,U:g_{\Upsilon_b})&\stackrel{\mathbf G_{bc}}{\rightarrow}&HM_k(f,U:g_{\Upsilon_c})&\rightarrow&\cdots\\
&&\downarrow \Phi^{\Upsilon^a} &&\downarrow \Phi^{\Upsilon^b}&&\downarrow \Phi^{\Upsilon^c}&&\\
\cdots&\rightarrow&HF_k(\Upsilon_a:H,J_a)&\stackrel{\mathbf F_{ab}}{\rightarrow}&HF_k(\Upsilon_b:H,J_b)&\stackrel{\mathbf F_{bc}}{\rightarrow}&HF_k(\Upsilon_c:H,J_c)&\rightarrow&\cdots\\
&&\uparrow \imath^{\lambda}_{\Upsilon^a*} &&\uparrow \imath^{\lambda}_{\Upsilon^b*}&&\uparrow \imath^{\lambda}_{\Upsilon^c*}&&\\
\cdots&\rightarrow&HF_k^\lambda(\Upsilon_a:H,J_a)&\stackrel{\mathbf F_{ab}^\lambda}{\rightarrow}&HF_k^\lambda(\Upsilon_b:H,J_b)&\stackrel{\mathbf F_{bc}^\lambda}{\rightarrow}&HF_k^\lambda(\Upsilon_c:H,J_c)&\rightarrow&\cdots
\end{array}
\end{equation}

Take $[\alpha]\in HM_k(f,U)\setminus\{0\}$ and $\lambda\in{\mathbb R}$ such that $\Phi([\alpha])\in\IM(\imath^\lambda_*)$; there exists $[x]\in HF_k^\lambda(H,U:M)$ such that
$$
\Phi([\alpha])=\imath^\lambda_*([x]).$$

From the definition of a direct limit we conclude that
$$
\alpha\in HM_k(f,U:g_\Upsilon),\hskip1mmx\in HF_k^\lambda(\Upsilon':H,J_{\Upsilon'})$$
for some $\Upsilon$ and $\Upsilon'$. Since
$$
\Phi([\alpha])=[\Phi^\Upsilon(\alpha)]=[\imath^\lambda_{\Upsilon'*}(x)]=\imath^\lambda_*[x],
$$
we find that
$$
F_{\Upsilon\overline{\Upsilon}}(\Phi^\Upsilon(\alpha))=F_{\Upsilon'\overline{\Upsilon}}(\imath^\lambda_{\Upsilon'*}(x))
$$
for some $\overline{\Upsilon}$ which is closer to $\nu^*\overline{U}$ than $\Upsilon$ and $\Upsilon'$, $\overline{\Upsilon}\leq\Upsilon$ and $\overline{\Upsilon}\leq\Upsilon'$. Using the commutativity~(\ref{eq:2diagrams}) we get
$$
\Phi^{\overline{\Upsilon}}(G_{\Upsilon\overline{\Upsilon}}(\alpha))=F_{\Upsilon\overline{\Upsilon}}(\Phi^\Upsilon(\alpha))=F_{\Upsilon'\overline{\Upsilon}}(\imath^\lambda_{\Upsilon'*}(x))
=\imath^{\lambda}_{\overline{\Upsilon}*}(F^\lambda_{\Upsilon'\overline{\Upsilon}}(x)).
$$
Therefore
$$
\Phi^{\overline{\Upsilon}}(G_{\Upsilon\overline{\Upsilon}}(\alpha))\in\IM(\imath^\lambda_{\overline{\Upsilon}*}).
$$
We conclude 
\begin{equation}\label{eq:spect1}
c_{\overline{\Upsilon}}(G_{\Upsilon\overline{\Upsilon}}(\alpha),H)\leq c_U([\alpha],H).
\end{equation}
If we take $\alpha\in HM_k(f,U:g_{\Upsilon})\setminus\{0\}$ and $\lambda\in{\mathbb R}$ such that $$\Phi^\Upsilon(\alpha)\in\IM(\imath^\lambda_{\Upsilon*}),$$ then $\Phi^\Upsilon(\alpha)=\imath^\lambda_{\Upsilon*}(x)$ for some $x\in HF_k^\lambda(\Upsilon:H,J_\Upsilon)$. Therefore, we have
$$
\Phi([\alpha])=[\Phi^\Upsilon(\alpha)]=[\imath^\lambda_{\Upsilon*}(x)]=\imath^\lambda_*[x],
$$
so we obtain the inequality
\begin{equation}\label{eq:spect2}
c_U([\alpha],H)\leq c_\Upsilon(\alpha,H).
\end{equation}
The elements $\alpha$ and $G_{\Upsilon\overline{\Upsilon}}(\alpha)$ represent the same element in the quotient space $HM_k(f,U)$. From (\ref{eq:spect1}) and (\ref{eq:spect2}) we have
\begin{equation}\label{eq:ineqaux}
c_{\overline{\Upsilon}}(G_{\Upsilon\overline{\Upsilon}}(\alpha),H)\leq c_U([\alpha],H)=c_U([G_{\Upsilon\overline{\Upsilon}}(\alpha)],H)\leq c_{\overline{\Upsilon}}(G_{\Upsilon\overline{\Upsilon}}(\alpha),H),
\end{equation}
so all inequalities become equalities.

Note that spectral invariants decrease as $\Upsilon\rightarrow\nu^*\overline{U}$, i.e.
\begin{equation}\label{eq:inequ_inv}
c_{\widetilde{\Upsilon}}(G_{\Upsilon\widetilde{\Upsilon}}(\alpha),H)\leq c_{\Upsilon}(\alpha,H),
\end{equation}
for every $\widetilde{\Upsilon}\le\Upsilon$. Indeed, if $\Phi^\Upsilon(\alpha)\in\IM(\imath^\lambda_{\Upsilon*})$, for $\alpha\neq 0\in HM_*(f,U;g_\Upsilon)$, we have
\begin{equation}\label{eq:auxili}
\Phi^\Upsilon(\alpha)=\imath^\lambda_{\Upsilon*}(x)\Rightarrow F_{\Upsilon\widetilde{\Upsilon}}(\Phi^\Upsilon(\alpha))=F_{\Upsilon\widetilde{\Upsilon}}(\imath^\lambda_{\Upsilon*}(x)),
\end{equation} so, from the commutativity of the diagram~(\ref{eq:2diagrams}) we have
$$\imath^\lambda_{\widetilde\Upsilon*}(F_{\Upsilon\widetilde{\Upsilon}}(x))=F_{\Upsilon\widetilde{\Upsilon}}(\imath^\lambda_{\Upsilon*}(x))\stackrel{(\ref{eq:auxili})}{=}F_{\Upsilon\widetilde{\Upsilon}}(\Phi^\Upsilon(\alpha))=\Phi^{\widetilde\Upsilon}(G_{\Upsilon\widetilde{\Upsilon}}(\alpha)).
$$
This means that 
$$\Phi^{\widetilde\Upsilon}(G_{\Upsilon\widetilde{\Upsilon}}(\alpha))\in\IM(\imath^\lambda_{\widetilde\Upsilon*}),$$ so~(\ref{eq:inequ_inv}) holds.

From~(\ref{eq:ineqaux}) and~(\ref{eq:inequ_inv}) we conclude that all $c_{\Upsilon}(\alpha,H)$ become equal to $c_U([\alpha],H)$, starting from some $\widetilde{\Upsilon}$.\qed

\bigskip

\noindent{\bf Proof of the part (E) in Theorem~\ref{thm:inv_open}.} In the same way as in the proof of  Lemma 2.6 in~\cite{MVZ} one can prove that the invariants for approximations are the same for two compactly supported Hamiltonians generating the same time-one-map, i.e:
$$\phi_H^1=\phi_K^1\,\Longrightarrow\, c_\Upsilon(\alpha,H)=c_\Upsilon(\alpha,K).$$
The proof now follows directly from Proposition~\ref{prop:limit}.\qed

\subsection{Continuity of spectral invariants}

The following theorem is the part (B) of Theorem~\ref{thm:inv_open}.
 
\begin{thrm}\label{thm:Cont_invar} Let $\|\cdot\|$ denotes the Hofer's norm:
$$\|H\|:=\int_0^1[\max_xH(x,t)-\min_xH(x,t)]dt.$$
Relative spectral invariants for an open set $U$
$$C_U(\alpha,H):=c_U(\alpha,H)-c_U(1,H)
$$
are continuous with respect to $\|\cdot\|:$
$$
|C_U(\alpha,H)-C_U(\alpha,H')|\le\|H-H'\|.
$$
(Here $1$ denotes the generator of zero homology group $HM_0(f,U:g_\Upsilon)$.)
 
\end{thrm}
\noindent{\it Proof:} First, we prove that $$C_\Upsilon(\alpha,H):=c_\Upsilon(\alpha,H)-c_\Upsilon(1,H)$$ is continuous with respect to Hamiltonian $H$. Let us fix a good approximation $\Upsilon$ and let $H$ and $H'$ be two Hamiltonians satisfying $(\ref{prop:trans_cond})$. Consider the linear homotopy $$H^s=(1-s)H+sH'=H+\sigma(s)(H'-H),$$ (we can approximate this linear homotopy with a regular one). The isomorphism $S^{\Upsilon}_{H,H'}$ is defined by a number of the holomorphic strips that connect $x\in CF_*(O_M,\Upsilon:H,J_\Upsilon)$ and $y\in CF_*(O_M,\Upsilon:H',J_\Upsilon)$:
$$\begin{aligned}&\mathcal{M}(x,y,O_M,\Upsilon:H,H',J_\Upsilon):=\\
&\left\{u:{\mathbb R}\times[0,1]\to T^*M\left|\begin{array}{l}
\frac{\partial u}{\partial s}+
J_\Upsilon(\frac{\partial u}{\partial t}-X_{H^s}(u))= 0\\
u(s,0)\in O_M,\,u(s,1)\in\Upsilon,  \\
u(-\infty,t)=x(t),\, u(+\infty,t)=y(t)\\
\end{array}\right.\right\}.\end{aligned}$$
If there exists $u\in\mathcal{M}(x,y,O_M,\Upsilon:H,H',J_\Upsilon)$, for a linear homotopy $H^s$, then by direct computation we see that it holds
\begin{equation}\label{eq:aux_est}
\mathcal{A}_{H'}^\Upsilon(y)-\mathcal{A}_{H}^\Upsilon(x)=\int_{-\infty}^{+\infty}\frac{d}{ds}\mathcal{A}_{H^s}^\Upsilon(u(s,\cdot))\leq E_+(H-H'):=\int_0^1\max_x(H-H')\,dt.
\end{equation}
Since linear homotopy may not be regular, we can approximate it by a $C^1$-close regular homotopy $H^s$, and obtain:
$$\mathcal{A}_{H'}^\Upsilon(y)-\mathcal{A}_{H}^\Upsilon(x)\leq E_+(H-H')+\varepsilon$$ for any $\varepsilon>0$. Letting $\varepsilon\to 0$, we get the estimate~(\ref{eq:aux_est}) for a regular homotopy $H^s$.
It follows
\begin{equation}\label{ineq_inv}\mathcal{A}_{H'}^\Upsilon(S_{H,H'}^\Upsilon(x))\leq\mathcal{A}_{H}^\Upsilon(x)+E_+(H-H').
\end{equation}
For $x\in HF_*(O_M,\Upsilon:H,J_\Upsilon)$ we can define
$$
\widetilde{c}_\Upsilon(x,H):=\inf\{\lambda\in{\mathbb R}\,|\,x\in \IM(\imath_{\Upsilon,H*}^\lambda)\}.
$$
Obviously, it holds:
$$c_\Upsilon(\alpha,H)=\widetilde{c}_\Upsilon(\Phi^\Upsilon_H(\alpha),H).$$
It follows from~(\ref{ineq_inv}):
$$
\widetilde{c}_\Upsilon(S_{H,H'}^\Upsilon(x),H')\leq\widetilde{c}_\Upsilon(x,H)+E_+(H-H').
$$
Since $S_{H,H'}^\Upsilon\circ\Phi^\Upsilon_H=\Phi^\Upsilon_{H'}$ we get the inequality
$$\begin{aligned}
c_\Upsilon(\alpha,H')&=\widetilde{c}_\Upsilon(\Phi^\Upsilon_{H'}(\alpha),H')=\widetilde{c}_\Upsilon(S_{H,H'}^\Upsilon\circ \Phi^\Upsilon_H(\alpha),H')\\
&\leq\widetilde{c}_\Upsilon(\Phi^\Upsilon_H(\alpha),H)+E_+(H-H')\\
&=c_\Upsilon(\alpha,H)+E_+(H-H'),
\end{aligned}$$
that holds for all $\alpha\in HM_*(f,U:g_\Upsilon)$. If we write the same inequality for the generator of zero homology group, we derive the continuity of relative spectral invariants for approximations:
$$|C_\Upsilon(\alpha,H')-C_\Upsilon(\alpha,H)|\leq\|H-H'\|.
$$

Now the proof follows from the above inequality and Proposition~\ref{prop:limit}.
\qed

\bigskip


\subsection{Triangle inequality}

Now we prove the part (A) in Theorem~\ref{thm:inv_open}.

For two function $H_1, H_2:T^*M\times[0,1]\to\mathbb{R}$ with $H_1(x,1)=H_2(x,0)$, we define their concatenation as:
$$H_1\sharp\,H_2:=\begin{cases}H_1(x,t),&t\le 1\\
H_2(x,t-1),&t\ge 1.
\end{cases}$$

\begin{prop}\label{prop:cont_app} Let $f_j\in\mathcal{F}^-(g)$, for $j=1,2,3$ (see Definition~\ref{defn:f}) and $\alpha\in HM_*(f_1,U:g)$, $\beta\in HM_*(f_2,U:g)$. If $0\neq\alpha\cdot\beta\in HM_*(f_3,U:g)$
then
$$
c^3_\Upsilon(\alpha\cdot\beta,H_1\sharp H_2)\le c^1_\Upsilon(\alpha,H_1)+c^2_\Upsilon(\beta,H_2)
$$ where $c^j_\Upsilon$ denotes the invariant defined via PSS isomorphism that involves Morse function $f_j$.
\end{prop}
\noindent{\it Proof:} Choose a Hamiltonian $H_3$ that is regular, smooth and close enough to $H_1\sharp H_2$:
$$\|H_3-H_1\sharp H_2\|_{C^0}<\varepsilon.
$$
We prove that a product $\circ$ descends to a product on filtered homologies
$$\circ:HF_*^\lambda(\Upsilon:H_1)\otimes HF_*^\sigma(\Upsilon:H_2)\longrightarrow HF_*^{\lambda+\sigma+4\varepsilon}(\Upsilon:H_3).
$$
Let $\Sigma$ denotes the Riemannian surface defined in Subsection~\ref{subsec:product}.
Take a smooth family of Hamiltonians $K:\Sigma\times T^*M\to{\mathbb R}$ such that
$$K(s,t,\cdot)=\begin{cases}H_1(t+1,\cdot),&s\le -1,-1\le t\le0\\
        H_2(t,\cdot),&s\le -1,0\le t\le1\\
        \frac{1}{2}H_3(\frac{t+1}{2},\cdot),&s\ge1.
        \end{cases}$$
We can choose $K$ such that
$$
\bigg\|\frac{\partial K}{\partial s}\bigg\|\le\varepsilon,\;s\in[-1,1],
$$
and
$$
\frac{\partial K}{\partial s}=0,
$$
elsewhere. Assume that for $x\in CF_*^\lambda(\Upsilon:H_1)$ and $y\in CF_*^\sigma(\Upsilon:H_2)$ there exists $u\in{\mathcal M}(x,y;z)$ for some $z\in CF_*(\Upsilon:H_3)$. Here, $u$ are pseudo--holomorphic pants for a Hamiltonian $K$
$$
\bar{\partial}_{K,J^{\Upsilon}}(u)=0.
$$
Using the relations
$$
\begin{aligned}
&\int_{\Sigma}\bigg\|\frac{\partial u}{\partial s}\bigg\|^2\,ds\,dt\ge0\\
&\int_\Sigma u^*\omega=-\int x^*\theta+h_{\Upsilon}(x(1))-\int y^*\theta+h_{\Upsilon}(y(1))+\int z^*\theta-h_{\Upsilon}(z(1)),
\end{aligned}
$$
Stoke's formula and properties of a Hamiltonian $K$ it follows
$$
{\mathcal A}^{\Upsilon}_{H_3}(z)\le{\mathcal A}^{\Upsilon}_{H_1}(x)+{\mathcal A}^{\Upsilon}_{H_2}(y)+4\varepsilon.
$$
Now, from Proposition~\ref{prop:circ_cdot_PSS} we obtain the inequality
$$
c^3_{\Upsilon}(\alpha\cdot\beta,H_3)\le c^1_\Upsilon(\alpha,H_1)+c^2_\Upsilon(\beta,H_2)+4\varepsilon.
$$
Since spectral invariants are continuous with respect to the Hamiltonian the proof follows.
\qed
\\

\begin{thrm} For $[\alpha],[\beta]\in HM_*(f,U)$ such that $[\alpha\cdot\beta]\neq0$ it holds
$$
c^3_U([\alpha\cdot\beta],H_1\sharp H_2)\le c^1_U([\alpha],H_1)+c^2_U([\beta],H_2),
$$
where $j$ in $c_U^j$ emphasizes the corresponding Morse function $f_j$.
\end{thrm}

\noindent{\it Proof.} From Proposition~\ref{prop:limit} we have
\begin{align*}
&c^1_U([\alpha],H)=c^1_{\overline{\Upsilon}}(G_{\Upsilon\overline{\Upsilon}}(\alpha),H)\\
&c^2_U([\beta],H)=c^2_{\overline{\Upsilon}}(G_{\Upsilon\overline{\Upsilon}}(\beta),H)\\
&c^3_U([\alpha\cdot\beta],H)=c^3_{\overline{\Upsilon}}(G_{\Upsilon\overline{\Upsilon}}(\alpha\cdot\beta),H)
\end{align*}
for all $\overline{\Upsilon}\le\widetilde{\Upsilon}$. For such an $\Upsilon$, it follows from Proposition~\ref{prop:cont_app}:
$$c^3_{\overline{\Upsilon}}(\mathbf{G}^{f_1}_{\Upsilon\overline{\Upsilon}}(\alpha)\cdot\mathbf{G}^{f_2}_{\Upsilon\overline{\Upsilon}}(\beta),H_3)\le c^1_{\overline{\Upsilon}}(\mathbf{G}^{f_1}_{\Upsilon\overline{\Upsilon}}(\alpha),H_1)+c^2_{\overline{\Upsilon}}(\mathbf{G}_{\Upsilon\overline{\Upsilon}}^{f_2}(\beta),H_2).$$ Now, from~(\ref{eq:morse_prod}) we have
$$\mathbf{G}_{\Upsilon\overline{\Upsilon}}^{f_3}(\alpha\cdot \beta)=\mathbf{G}_{\Upsilon\overline{\Upsilon}}^{f_1}(\alpha)\cdot\mathbf{G}_{\Upsilon\overline{\Upsilon}}^{f_2}(\beta)$$ and, therefore
$$\begin{aligned}
c^3_U([\alpha\cdot\beta],H_3)&=c^3_{\overline{\Upsilon}}(\mathbf{G}^{f_3}_{\Upsilon\overline{\Upsilon}}(\alpha\cdot\beta),H_3)\le 
c^1_{\overline{\Upsilon}}(\mathbf{G}^{f_1}_{\Upsilon\overline{\Upsilon}}(\alpha),H_1)+c^2_{\overline{\Upsilon}}(\mathbf{G}^{f_2}_{\Upsilon\overline{\Upsilon}}(\beta),H_2)
\\
&=c^1_U([\alpha],H_1)+c^2_U([\alpha],H_2).\end{aligned}$$
\qed

\bigskip

\subsection{Invariants for periodic orbits}\label{subsec:per}
Recall the definition of spectral invariants for periodic orbit Floer homology. Since Floer homology for periodic orbits is not well defined for compactly supported Hamiltonians in $T^*M$, we will consider Hamiltonians with a support in some fixed cotangent ball bundle as in~\cite{FS} and also used in~\cite{MVZ}. More precisely, fix $R>0$, $\varepsilon>0$ and a smooth function $h:(-\varepsilon,+\infty)\to\mathbb{R}$ with the following properties:
\begin{itemize}
\item $h(t)=0$ for $t\ge 0$;
\item $h'(t)\ge 0$ for $t\le 0$;
\item $h'$ is small enough so that the flow of $h(\|p\|-R)$ does not have non constant periodic orbit of period less or equal to $1$ for $\|p\|\in(0,\varepsilon)$.
\end{itemize}
We choose $H_t(q,p)$ to be equal to $h(\|p\|-R)$ for $\|p\|\ge R-\varepsilon$.

Denote by $HF_k(T^*M:H,J)$ and $HM_k(f,T^*M)$ Floer homology for periodic orbits in $T^*M$ and Morse homology for the Morse function $f:T^*M\to\mathbb{R}$ respectively. Denote by $HF_k^{\lambda}(T^*M:H,J)$ the corresponding filtered group (with respect to the standard action functional) and, again, by $\jmath_*^\lambda$ the map induced by the inclusion map. Let $\PSS$ stands for PSS isomorphism for periodic orbits, defined in a way analogous to~\cite{PSS}
$$\PSS:HM_k(f,T^*M:g)\stackrel{\cong}{\longrightarrow}HF_k(T^*M:H,J)$$ and let $\alpha\in HM_k(f,T^*M:g)$.

The filtration in Floer homology for periodic orbits is given by the standard action functional
$$a_H(\gamma):=\int\gamma^*\theta-\int_0^1 Hdt$$
which is well defined in the cotangent bundle setting. Filtered Floer homology groups are homology groups of a chain complex generated by
$$CF_k^{\lambda}(T^*M:H):=\{a\in CF_k(T^*M:H)\mid a_H(a)<\lambda\},$$
where $CF_k(T^*M:H)$ denotes the $\mathbb Z_2-$vector space over the set of periodic Hamiltonian $H-$orbits in $T^*M$ of Conley--Zehnder index $k$.
\begin{defn} Let $\alpha\in HM_*(f,T^*M:g)\setminus\{0\}$.
Define
$$\rho(\alpha,H):=
\inf\{\lambda\mid\PSS(\alpha)\in\IM(\jmath^\lambda_*)\}.$$
\end{defn}

\subsection{Chimneys and relation between the two invariants}\label{subsec:chim}
This subsection is dedicated to a comparison of spectral invariants in periodic orbits and Lagrangian case and the proof of the part (C) of Theorem~\ref{thm:inv_open}.

The homomorphisms defined using "chimneys" are constructed by Abbondandolo and Schwarz in~\cite{AS} (in the context of Floer homology of cotangent bundles and the ring-isomorphism with the homology of the loop space) and Albers in~\cite{A} (in the construction of the comparison homomorphisms between Lagrangian and Hamiltonian Floer homology). The construction of a chimney is different in our situation, due to the boundary conditions.

Let
$$\Pi:=\mathbb{R}\times [0,1]/\sim,\quad\mbox{where}\; (s,0)\sim (s,1)\;\mbox{for}\;s\ge 0.$$

\begin{figure}
\centering
\includegraphics[width=8cm,height=4cm]{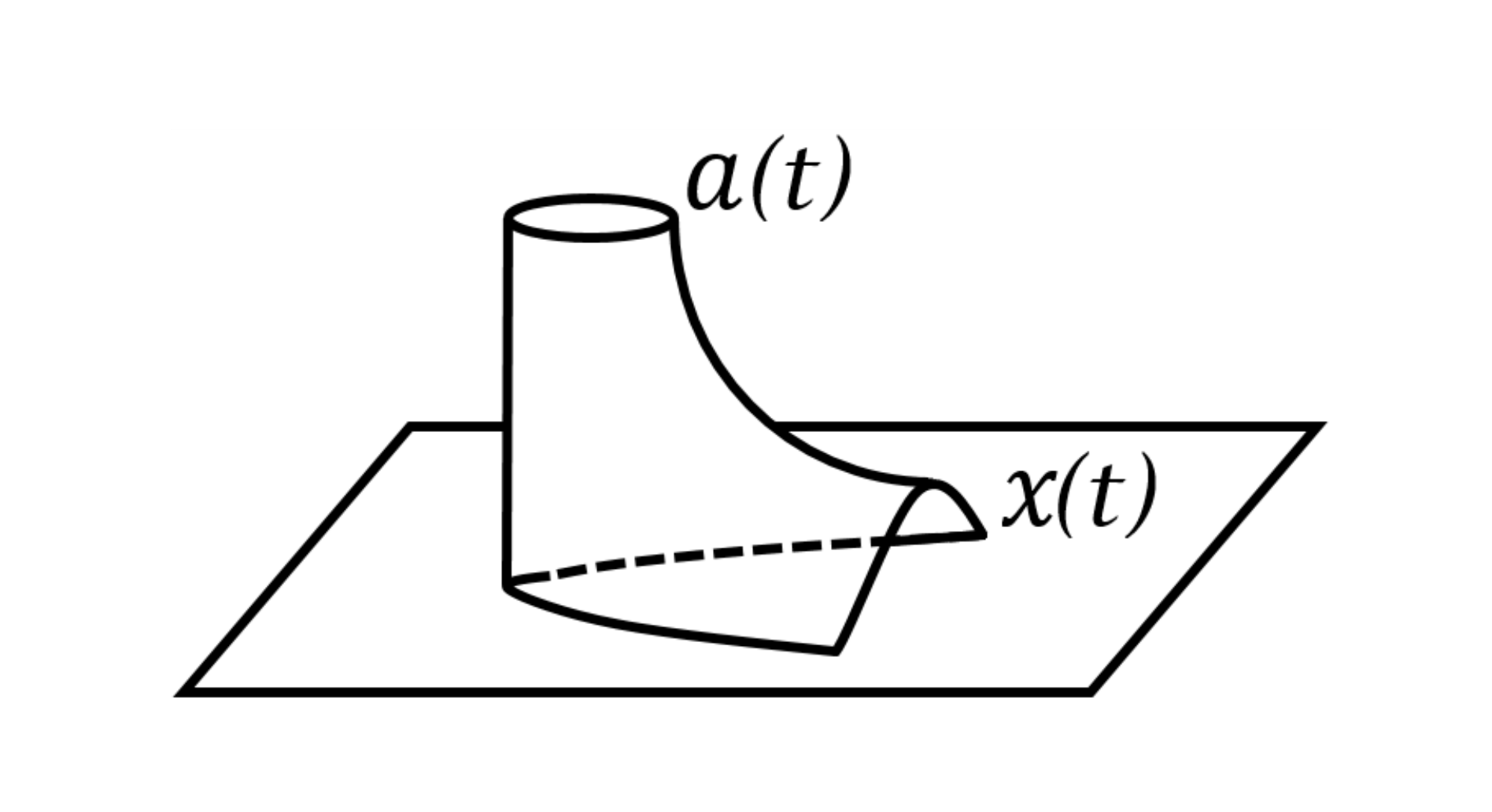}
\centering
\caption{Chimney}
\label{fig:chimney}
\end{figure}

For $x\in CF_*(O_M,\Upsilon:H)$ and $a\in CF_*(T^*M:H)$ define the manifold
of chimneys as:
$$
\mathcal{M}(x,a,O_M,\Upsilon:H,J):=
\left\{u:\Pi\to T^*M\left|\begin{array}{l}\partial_su+J(\partial_tu-X_H\circ u)=0\\
u(s,0)\in O_M,\,u(s,1)\in\Upsilon\;\mbox{for}\;s\le 0\\
u(-\infty,t)=x(t), \,u(+\infty,t)=a(t)
\end{array}\right.\right\}
$$
(see Figure~\ref{fig:chimney}). For generic choices, $\mathcal{M}(x,a)$ is a smooth manifold of dimension $\mu_{CZ}(a)-\mu(x)-\frac{n}{2}$, where $\mu_{CZ}(a)$ denotes the Conley-Zehnder index of a loop $a$.

Define
\begin{equation}\label{eq:chi}
\begin{aligned}&\chi:CF_k(O_M,\Upsilon:H)\to CF_k(T^*M:H)\\
&\chi(x):=\sum\limits\sharp_2\,\mathcal{M}(x,a,O_M,\Upsilon:H,J)\,a.
\end{aligned}\end{equation}
It holds $\chi\circ\partial=\partial\circ\chi$, hence $\chi$ is well defined on the homology level:
$$\chi:HF_k(O_M,\Upsilon:H,J)\to HF_k(T^*M:H,J).$$

Let $a$ be a periodic orbit. If there exists $u\in\mathcal{M}(x,a,O_M,\Upsilon:H,J)$, let $y$ be a periodic orbit defined as
$$y(t):=u(0,t).$$ Since $y(0)=y(1)\in O_M$, we have $h_{\Upsilon}(y(1))=0$. Therefore, we have
$$\mathcal{A}_H^\Upsilon(y)=a_H(y),$$
so
$$\begin{aligned}
&a_H(a)-\mathcal{A}_H^\Upsilon(x)=a_H(a)-a_H(y)+\mathcal{A}_H^\Upsilon(y)-\mathcal{A}_H^\Upsilon(x)=\\
&\int_{-\infty}^{0}\frac{d}{ds}a_H(u(s,t))ds+
\int_{0}^{+\infty}\frac{d}{ds}\mathcal{A}_H^\Upsilon(u(s,t))ds=\\
&-\int_{-\infty}^{\infty}\int_0^1\omega\left(\partial_su,\partial_tu-X_H\circ u\right)\,dtds
=-\int_{-\infty}^{\infty}\int_0^1\left\|\frac{\partial u}{\partial s}\right\|^2\,dtds
\le 0.\end{aligned}$$
It follows that $\chi$ defines the mapping
$$\chi^\lambda:=\chi|_{CF^{\lambda}_k(O_M,\Upsilon:H)}:
CF^{\lambda}_k(O_M,\Upsilon:H)\to CF^{\lambda}_k(T^*M:H)$$
which also descends to the homology level:
$$\chi^\lambda:
HF^{\lambda}_k(O_M,\Upsilon:H,J_\Upsilon)\to HF^{\lambda}_k(T^*M:H,J_\Upsilon)$$
(see also~\cite{DKM}).
The diagram
\begin{equation}\label{eq:diag:chi,i,j}\xymatrix{
 HF^{\lambda}_k(O_M,\Upsilon:H,J_\Upsilon)\ar[d]_{\imath_*^{\lambda}}\ar[r]^-{\chi^{\lambda}} &HF^{\lambda}_k(T^*M:H,J_\Upsilon)
  \ar[d]_{\jmath_*^{\lambda}}\\
 HF_k(O_M,\Upsilon:H,J_\Upsilon)\ar[r]^-{\chi} &HF_k(T^*M:H,J_\Upsilon)}\end{equation} commutes.

Similarly, set
$$\Delta=\mathbb{R}\times [0,1]/\sim,\quad\mbox{where}\; (s,0)\sim (s,1)\;\mbox{for}\;s\le 0,$$
and define

$$
\mathcal{M}(a,x,O_M,\Upsilon:H,J):=
\left\{u:\Delta\to T^*M\left|\begin{array}{l}\partial_su+J(\partial_tu-X_H\circ u)=0\\
u(s,0)\in O_M,\,u(s,1)\in\Upsilon\;\mbox{for}\;s\ge 0\\
u(-\infty,t)=a(t), \,u(+\infty,t)=x(t).
\end{array}\right.\right\}
$$

For generic choices, $\mathcal{M}(a,x)$ is a smooth manifold of dimension $\mu(x)-\frac{n}{2}-\mu_{CZ}(a)$.

Define
$$\begin{aligned}&\tau:CF_k(T^*M:H)\to CF_{k-n}(O_M,\Upsilon:H)\\
&\tau(a):=\sum\limits\sharp_2\,\mathcal{M}(a,x,O_M,\Upsilon:H,J)\,x.
\end{aligned}$$ This homomorphism also
descends to the homology level
$$\tau:HF_k(T^*M:H,J)\to HF_{k-n}(O_M,\Upsilon:H,J)$$ since it commutes with the boundary operators. As above, one can show that it also induces a homomorphism
on the filtered homology level, and that the corresponding diagram (analogous to~(\ref{eq:diag:chi,i,j})) commutes.

Let $f\in\mathcal{F}^-(g)$ (see Definition~\ref{defn:f} on page~\pageref{F^+}). We extend $f$ to $F:T^*M\to\mathbb{R}$ in the following way. Consider a tubular neighbourhood $W \subseteq T^*M$ of $M$. First, we extend $f$ to the vector bundle $W$ over $M$, and obtain the Morse function $f_W$ such that
$$ f_W|_M=f\,\text{and}\,\,\,\Crit(f_W)=\Crit(f).
$$
Then we extend the Morse function $f_W$ defined on the open subset $W$ to the Morse function $F$ on $T^*M$ such that there are no trajectories for the negative gradient flow of $F$ leaving $W$ (see \cite{Sc1} for details). Now the Morse complex $CM_*(f)$ is a subset of the Morse complex $CM_*(F)$ and the inclusion map of these complexes becomes the homomorphism $\imath_*$ on the homology level.

\begin{prop} Let $f$, $F$ and $\imath_*$ be as above. Let $H:T^*M\to\mathbb{R}$ be a Hamiltonian. Suppose all the choices are generic. The diagram
\begin{equation}\label{2diags}
\xymatrix{ HF_k^{\lambda}(O_M,\Upsilon:H,J_\Upsilon)\ar[d]_{\imath_*^{\lambda}}\ar[r]^-{\chi^{\lambda}} &HF_k^{\lambda}(T^*M:H,J_\Upsilon)
\ar[d]_{\jmath_*^{\lambda}}\\
HF_k(O_M,\Upsilon:H,J_\Upsilon)\ar[r]^-{\chi} &HF_k(T^*M:H,J_\Upsilon)
\ar[d]_{\PSS^{-1}}\\
HM_k(f,U:g_\Upsilon)\ar[u]_{\Phi^\Upsilon}\ar[r]^{\imath_*} &HM_k(F,T^*M:g_\Upsilon)}
\end{equation}
commutes.
\end{prop}

\noindent{\it Proof:}
The upper diagram is~(\ref{eq:diag:chi,i,j}). The lower diagram is
\begin{equation}\label{side_by_side}
\xymatrix{HF_k(O_M,\Upsilon:H,J_\Upsilon)\ar[r]^-{\chi} &HF_k(T^*M:H,J_\Upsilon)\ar[d]_{\PSS^{-1}}\\
HM_k(f,U:g_\Upsilon)\ar[u]_{\Phi^\Upsilon}\ar[r]^{\imath_*} &HM_k(F,T^*M:g_\Upsilon)}
\end{equation}
and its commutativity means that is holds
$$\PSS^{-1}\circ\chi\circ\Phi^\Upsilon=\imath_*.$$
In order to do that using the usual cobordism arguments, we consider the following two auxiliary manifolds. 

\begin{figure}
\centering
\includegraphics[width=9cm,height=5.5cm]{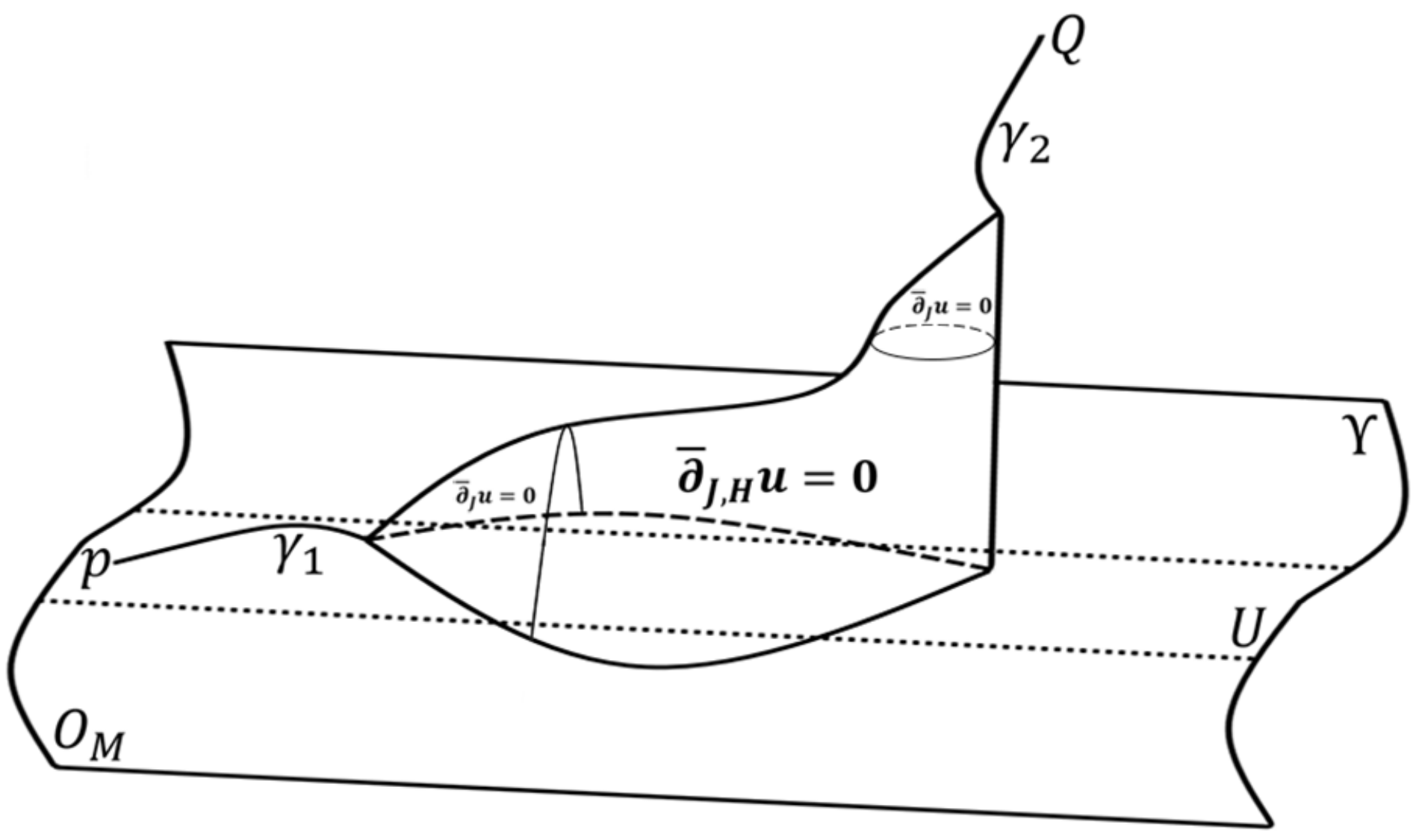}
\centering
\caption{Moduli space $\widetilde{\mathcal{M}}_R(p,Q)$}
\label{Mixed_pic}
\end{figure}

For $R>1$ let $\rho_R$ be a smooth cut-off function with the properties:
$$\rho_R(t)=\begin{cases} 0, & t\in(-\infty,-2R-1]\cup [-R+1,R-1]\cup [2R+1,+\infty)
 \\ 1, & t\in[-2R,-R]\cup[R,2R]. \end{cases}$$
Let $p$ be the critical point of a Morse function $f$ and $Q$ the critical point of a Morse function $F$. Fix $R>1$ and define
$$\begin{aligned}&
\widetilde{\mathcal{M}}_R(p,Q):=\widetilde{\mathcal{M}}_R(p,Q,O_M,\Upsilon:H,J_\Upsilon):=\\
&\left\{(\gamma_1,u,\gamma_2)\left|
\begin{array}{l}
\gamma_1:(-\infty,0]\to U,\;\dot{\gamma_1}=-\nabla f(\gamma_1)\\
\gamma_2:[0,+\infty)\to T^*M,\;\dot{\gamma_2}=-\nabla F(\gamma_2)\\
u:\Pi\to T^*M,\;\partial_su+J(\partial_tu-X_{\rho_R(s)H}\circ u)=0\\
u(s,0)\in O_M,\,u(s,1)\in\Upsilon\;\mbox{for}\;s\le 0\\
\gamma_1(-\infty)=p,\,
\gamma_1(0)=u(-\infty,t)\\
u(+\infty,t)=\gamma_2(0),\,
\gamma_2(+\infty)=Q
\end{array}\right.\right\}\end{aligned}$$
(see Figure~\ref{Mixed_pic}).

\noindent Define also
$$\widetilde{\mathcal{M}}(p,Q):=\left\{(\gamma_1,u,\gamma_2,R)\mid R\in [R_0,\infty),\;(\gamma_1,u,\gamma_2)\in
\widetilde{\mathcal{M}}_R(p,Q)\right\}.$$
For $m_f(p)=m_F(Q)$ and generic choices, $\widetilde{\mathcal{M}}(p,Q)$ is a smooth one-dimensional manifold with topological boundary that can be identified with
$$\partial\left(\widetilde{\mathcal{M}}(p,Q)\right)=
\mathcal{B}_1\cup \mathcal{B}_2\cup\mathcal{B}_3\cup\mathcal{B}_4$$
where
$$\begin{aligned}
&\mathcal{B}_1=\widetilde{\mathcal{M}}_{R_0}(p,Q);\\
&\mathcal{B}_2=\bigcup_{s\in\Crit(f)}\widehat{\mathcal{M}}(p,s)\times\widetilde{\mathcal{M}}_R(s,Q);\\
&\mathcal{B}_3=\bigcup_{S\in\Crit(F)}\widetilde{\mathcal{M}}_R(p,S)\times\widehat{\mathcal{M}}(S,Q);\\
&\mathcal{B}_4=\bigcup_{\scriptsize{\begin{array}{c}
x\in CF_k(O_M,\Upsilon:H)\\
a\in CF_k(T^*M:H)
\end{array}}}
\mathcal{M}(p,x)\times\mathcal{M}(x,a)\times\mathcal{M}(a,Q).
\end{aligned}$$
Here
$$\mathcal{M}(a,Q):=\left\{(u,\gamma)\left|
\begin{array}{l}
u:(-\infty,0]\times [0,1]\to T^*M\\
\partial_su+J(\partial_tu-X_{\tilde\rho H}\circ u)=0\\
u(s,0)=u(s,1)\\
\gamma:[0+\infty)\to T^*M,\,\dot{\gamma}=-\nabla F(\gamma)\\
u(0,t)=\gamma(0),\, u(-\infty,t)=a(t),\, \gamma(+\infty)=Q
\end{array}\right\},
\right.$$ i.e. it is the space of combined object defining a $\PSS$ isomorphism for periodic orbits and  $\tilde\rho(s)=\begin{cases}1,&s\le -2R-1\\ 0, &s\ge -2R\end{cases}$.

The boundary components $\mathcal{B}_2$ and $\mathcal{B}_3$ correspond to the boundary of $\widetilde{\mathcal{M}}_R(p,Q)$ (since $\widetilde{\mathcal{M}}(p,Q)\subset\widetilde{\mathcal{M}}_R(p,Q)\times[R_0,+\infty)$). The boundary parts $\mathcal{B}_1$ and $\mathcal{B}_4$ come from the $R$ coordinate, and $\mathcal{B}_1$ corresponds to the case $R\to R_0$. As regards the component $\mathcal{B}_4$, it arises when $R\to+\infty$. More precisely, for
$$
(\gamma_1,u_1)\in\mathcal{M}(p,x),\quad
v\in\mathcal{M}(x,a),\quad
(u_2,\gamma_2)\in\mathcal{M}(a,Q)$$
we define (for $R\ge 2$): 
\begin{equation}\label{eq:gluing}
u_1\,\sharp_R\,v\,\sharp_R\,u_2:=\left\{
\begin{array}{ll}
u_1(s+3R,t),&s\le-2R-1\\
\exp_{x(t)}(\beta(-s-2R)\xi_1(s+3R,t)),&-2R-1\le s\le -2R\\
x(t),&-2R\le s\le -R\\
\exp_{x(t)}(\beta(s+R)\eta_1(s-R,t)),& -R\le s\le -R/2\\
u(s,t),&-R/2\le s\le R/2\\
\exp_{a(t)}(\beta(-s+R)\eta_2(s+R)),&R/2\le s\le R\\
a(t),&R\le s\le 2R\\
\exp_{a(t)}\beta(s-2R)\xi_2(s-3R,t)),&2R\le s\le 2R+1\\ 
u_2(s-3R,t),&s\ge 2R+1
\end{array}\right.\end{equation}
and the approximative solution from Floer's gluing construction to be:
$$ (\gamma_1,u_1)\,\sharp_R\,
v\,\sharp_R\,(u_2,\gamma_2):=(\gamma_1,u_1\,\sharp_R\,v\,\sharp_R\,u_2,\gamma_2).$$
In the equation~(\ref{eq:gluing}) $\beta:\mathbb{R}\to[0,1]$ is a smooth cut-off function equal to $1$ for $s\ge 1$ and to $0$ for $s\le 0$. Vector fields $\xi_1(s,t),\eta_1(s,t)\in T_{x(t)}T^*M$ are chosen such that $$\begin{array}{cccc}u_1(s,t)&=\exp_{x(t)}, &t\in[0,1], &s\ge s_0\\
v(s,t)&=\exp_{x(t)}, &t\in[0,1], &s\le -s_0
\end{array}$$ and $\xi_2(s,t),\eta_2(s,t)\in T_{a(t)}T^*M$ are chosen similarly. The rest of the proof of Floer gluing theorem is standard.

Denote by $F_j$ the homomorphism obtained by counting the elements of the zero dimensional manifold $\mathcal{B}_j$. Since the number of the boundary of the one-dimensional manifold $\widetilde{\mathcal{M}}(p,Q)$ is even, i.e. zero in $\mathbb{Z}_2$, and the maps $F_2$ and $F_3$ are of the form
$$F_2=\partial\circ K,\quad F_3=K\circ\partial,$$
the homomorphisms $F_1$ and $F_4$ are equal in the homology. By standard cobordism argument one can show that the mapping $F_1$ does not depend on $R_0$. Now as in the proof of Theorem~\ref{thm:iso_approx}, we conclude that $F_1$ is chain homotopic to the map defined by the number of pairs $(\gamma_1,\gamma_2)$ with properties:
$$\left\{\begin{array}{lll}
\gamma_1:(-\infty,0]\to U,&\dot{\gamma_1}=-\nabla f(\gamma_1),&\gamma_1(-\infty)=p\\
\gamma_2:[0,+\infty)\to T^*M,
&\dot{\gamma_2}=-\nabla F(\gamma_2),&\gamma_2(+\infty)=Q\\
\gamma_1(0)=\gamma_2(0).&&
\end{array} \right.$$ Since $F|_U=f$, $\gamma_1\sharp\gamma_2$ is a negative gradient trajectory of $F$ connecting two critical points of the same Morse index. Thus, $F_1$ is chain homotopic to the homomorphism $\imath_0$. On the other hand, the mapping $F_4$ is exactly the homomorphism $\PSS^{-1}\circ\chi\circ\Phi^\Upsilon$, so the claim follows.

\qed
\bigskip

We intend to compare spectral invariants for two Floer homologies. Since in Lagrangian case we are dealing with the direct limit construction, i.e. we have the whole family of Floer homology groups (for the approximations) to start with, we need to have the corresponding family in periodic orbits case, to maintain the transversality conditions. In periodic orbit case, the canonical isomorphisms for two different almost complex structures will be the homomorphisms that define a direct limit Floer homology group.

For two generic almost complex structures $J_a$ and $J_b$, denote by $\mathbf{D}_{ab}$ a canonical isomorphism of Floer homologies for periodic orbits:
$$\mathbf{D}_{ab}:HF_k(T^*M:H,J_a)\to HF_k(T^*M:H,J_b)$$
that satisfies
$$\mathbf{D}_{bc}\circ \mathbf{D}_{ab}=\mathbf{D}_{ac}.$$ As before, define Floer homology for periodic orbits as a direct limit
$$HF_k(T^*M:H):=\limarr HF_k(T^*M:H,J_s).$$
The filtered Floer homology is defined as:
$$HF_k^{\lambda}(T^*M:H):=\limarr HF^{\lambda}_k(T^*M:H,J_s).$$

\begin{figure}
\centering
\includegraphics[width=9cm,height=3.6cm]{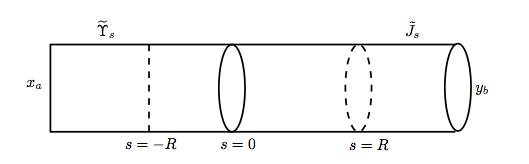}
\centering
\caption{Chimney from auxiliary manifold from Proposition~\ref{prop:com_chim_dir.lim}}
\label{aux_chimney_pic}
\end{figure}

\begin{prop}\label{prop:com_chim_dir.lim} Let $\chi^a$ stands for a homomorphism~(\ref{eq:chi}) for the almost complex structure $J_a$. We use the abbreviations
$$HF^{\lambda}_k(\Upsilon_a):=HF^{\lambda}_k(O_M,\Upsilon_a:H,J_a),\quad
HF^{\lambda}_k(J_a):=HF^{\lambda}_k(T^*M:H,J_a).$$
The diagram:
\begin{equation}\label{eq:diag_chi_ab}
\begin{array}{lllllllll}
\cdots&\longrightarrow&HF_k^{\lambda}(\Upsilon_a)&\stackrel{\mathbf F_{ab}}{\longrightarrow}&HF_k^{\lambda}(\Upsilon_b)&\stackrel{\mathbf F_{bc}}{\longrightarrow}&HF_k^{\lambda}(\Upsilon_c)&\longrightarrow&\cdots\\
&&\downarrow\chi^a&&\downarrow\chi^b&&\downarrow\chi^c&&\\
\cdots&\longrightarrow&HF^{\lambda}_k(J_a)&\stackrel{\mathbf D_{ab}}{\longrightarrow}&HF^{\lambda}_k(J_b)&\stackrel{\mathbf D_{bc}}{\longrightarrow}&HF^{\lambda}_k(J_c)&\longrightarrow&\cdots
\end{array}
\end{equation}
commutes.
\end{prop}
\noindent{\it Proof.} The commutativity of~(\ref{eq:diag_chi_ab}) is equivalent to:
$$\chi^a=\mathbf{D}_{ab}^{-1}\circ\chi^b\circ\mathbf{F}_{ab}=\mathbf{D}_{ba}\circ\chi^b\circ\mathbf{F}_{ab}.$$
The proof of the above equality is similar to the proofs of the Proposition~\ref{prop:var} and Lemma~\ref{lem}. The auxiliary one-dimensional manifold will be the set of the pairs $(R,u)$, where $R\in[R_0,+\infty)$, and $u$ is a chimney with the properties depicted in Figure~\ref{aux_chimney_pic}.
\qed

\begin{cor}
The homomorphism~(\ref{eq:chi}) induces the homomorphism
$$\chi^{\lambda}:HF_k^{\lambda}(H,U:M)\to HF_k^{\lambda}(T^*M:H),$$ and the homomorphism
$$\chi:HF_k(H,U:M)\to HF_k(T^*M:H).$$\qed
\end{cor}

The following corollary follows from the commutativity of~(\ref{2diags}) for all the approximations.

\begin{cor}\label{cor:com_dir_lim} The diagram
\begin{equation}\label{2diags_U}
\xymatrix{ HF_k^{\lambda}(H,U:M)\ar[d]_{\imath_*^{\lambda}}\ar[r]^-{\chi^{\lambda}} &HF_k^{\lambda}(T^*M:H)
\ar[d]_{\jmath_*^{\lambda}}\\
HF_k(H,U:M)\ar[r]^-{\chi} &HF_k(T^*M:H)
\ar[d]_{\PSS^{-1}}\\
HM_k(f,U)\ar[u]_{\Phi}\ar[r]^{\imath_*} &HM_k(F,T^*M)}
\end{equation} commutes.\qed
\end{cor}

\begin{thrm}\label{thm:comp1} Let $\alpha\in HM_k(f,U)\setminus\{0\}$. Then
$$c_U(\alpha,H)\ge\rho(\imath_*(\alpha),H).$$
\end{thrm}

\noindent{\it Proof:} From the commutativity of~(\ref{2diags_U})
ones easily gets
$$\left\{\lambda\mid \Phi(\alpha)\in\IM(\imath_*^{\lambda})\right\}
\subseteq
\left\{\lambda\mid \PSS(\imath_*(\alpha))\in\IM(\jmath_*^{\lambda})\right\},$$ so the claim follows.\qed
\bigskip

One can obtain the inequality of similar type by using the homomorphism $\tau$. The corresponding commutative diagram is
$$\xymatrix{ HF_k^{\lambda}(T^*M:H)\ar[d]_{\jmath_*^{\lambda}}\ar[r]^-{\tau^{\lambda}} &HF_{k-n}^{\lambda}(H,U:M)
\ar[d]_{\imath_*^{\lambda}}\\
HF_k(T^*M:H)\ar[d]_{\PSS^{-1}}\ar[r]^-{\tau} &HF_{k-n}(H,U:M)
\\
HM_k(F,T^*M)\ar[r]^{\imath_!} &HM_{k-n}(f,U)\ar[u]_{\Phi}},$$
where $\imath_!$ is the map obtained by inclusion map and Poincar\'e duality map:
$$\imath_!:=\PD^{-1}\circ\; \imath_*\circ\PD.$$
From this commutativity, we have the following
\begin{thrm} Let $\alpha\in HM_k(F,T^*M)\setminus\{0\}$, then
$$\rho(\alpha,H)\ge c_U(\imath_!(\alpha),H).$$
\end{thrm}
Proof is analogous to this of Theorem~\ref{thm:comp1}.\qed

\subsection{A remark on invariants for subsets}

In~\cite{O3} Oh considered a spectral invariant
$$c_+(H,U):=\inf\{\lambda\in\mathbb{R}\mid \imath^\lambda_*:HF^\lambda_k(H,U:M)\to HF_k(H,U:M)\;\mbox{is surjective}\}$$
(the notions are the same as in the Subsection~\ref{subsec:inv_open}). If $U\stackrel{\imath}{\hookrightarrow}V$ are two open subset of $M$ and
$$\imath_{*UV}:H^{\sing}_k(U,\mathbb{Z})\to H^{\sing}_k(V,\mathbb{Z})$$ is surjective, Oh proved that
$$c_+(H,V)\le c_+(H,U).$$
We can prove slightly more precise statement, the inequality for any homology class (with $\mathbb{Z}_2$ coefficients), using the PPS isomorphism for an open subset in the proof.

\begin{thrm}\label{thm:c_UV} Let $U\stackrel{\imath}{\hookrightarrow}V$ be two open subset of $M$ and let
$$\jmath_{*UV}:HM_k(f,U)\to HM_k(f,V)$$ (the homomorphism induced by inclusion $\jmath:U\hookrightarrow V$) be surjective. Let $c_U(\alpha,H)$ be as in~(\ref{eq:def_inv_c_open}). For $\alpha\in HM_k(f,U)\setminus\{0\}$ it holds:
$$c_V(\jmath_{*UV}(\alpha),H)\le c_U(\alpha,H).$$
\end{thrm}

\noindent{\it Proof:} Let
$$\begin{aligned}
&\imath_{*UV}:HF_k(H,U:M)\to HF_k(H,V:M)\\
&\imath^\lambda_{*UV}:HF^\lambda_k(H,U:M)\to HF^\lambda_k(H,V:M)\\
&\imath^\lambda_{*U}:HF^\lambda_k(H,U:M)\to HF_k(H,U:M)\\
&\imath^\lambda_{*V}:HF^\lambda_k(H,V:M)\to HF_k(H,V:M)
\end{aligned}$$ denote the inclusion homomorphisms defined by Oh in~\cite{O3}.
The following diagram is commutative:
\begin{equation}\label{eq:2diags}
\xymatrix{ HF_k^{\lambda}(H,U:M)\ar[d]_{\imath_{*U}^{\lambda}}\ar[r]^-
{\imath_{*UV}^{\lambda}}
&HF_k^{\lambda}(H,V:M)
\ar[d]_{\imath_{*V}^{\lambda}}\\
HF_k(H,U:M)\ar[d]_{\Psi_U}\ar[r]^-{\imath_{*UV}} &HF_k(H,V:M)
\ar[d]_{\Psi_V}\\
HM_k(f,U)\ar[r]^{\jmath_{*UV}} &HM_k(f,V).}
\end{equation}
The commutativity of the upper diagram is proven in~\cite{O3}. To prove the commutativity of the lower one, it is enough to prove the commutativity of
$$\xymatrix{
HF_k(O_M,\Upsilon^U:H,J)\ar[d]_{\Psi_{\Upsilon^U}}\ar[r]^-{\imath^{(H,J)}_{*UV}} &HF_k(O_M,\Upsilon^V:H,J)
\ar[d]_{\Psi_{\Upsilon^V}}\\
HM_k(f,U:g)\ar[r]^{\jmath_{*UV}} &HM_k(f,V:g),}$$ for all $\Upsilon^U$ close enough to $\nu^*\overline{U}$ and  $\Upsilon^V$ close enough to $\nu^*\overline{V}$. Here $\imath^{(H,J)}_{*UV}$ is the inclusion map also defined in~\cite{O3}. Take $[x]$ in $HF_k(O_M,\Upsilon^U:H,J)$. It holds
\begin{equation}\label{eq:aux4}
\Psi_{\Upsilon^V}(\imath^{(H,J)}_{*UV}([x]))=\sum_{p\in CM_k(V)}n(x,p)[p].\end{equation} On the other hand, we have
$$\jmath_{*UV}(\Psi_{\Upsilon^U}([x]))=
\sum_{p\in CM_k(U)}n(x,p)\jmath_{*UV}([p]),$$
which is the same as~(\ref{eq:aux4}) if $\jmath_{*UV}$ is surjective.

Let
$$A^U_\alpha:=\{\lambda\in\mathbb{R}\mid \Phi_U(\alpha)\in\IM(\imath^\lambda_{*U})\}.$$ If $\lambda\in A^U_\alpha$, then
$\Phi_U(\alpha)=\imath^\lambda_{*U}(\beta)$, for $\beta\in HF^\lambda_k(H,U:M)$, so, from the commutativity of~(\ref{eq:2diags}) we have
$$\imath^\lambda_{*V}(\imath^\lambda_{*UV}(\beta))=
\imath_{*UV}(\imath^\lambda_{*U}(\beta))=\imath_{*UV}(\Phi_U(\alpha))=\Phi_V(\jmath_{*UV}(\alpha)).$$
We conclude that $\lambda\in A^V_{\jmath_{*UV}(\alpha)}$, therefore
$$A^U_\alpha\subset A^V_{\jmath_{*UV}(\alpha)},$$ so by taking an infimum over $\lambda$, we obtain
$$c_V(\jmath_{*UV}(\alpha),H)\le c_U(\alpha,H).$$
\qed

\end{document}